\documentclass[11pt]{article}
\usepackage{amsmath, graphicx, amsfonts,amssymb}
\usepackage{amsfonts,mathrsfs, color, amsthm}
\usepackage{bm}

\usepackage{enumitem}
\usepackage{hyperref}
\hypersetup{
	colorlinks   = true, 
	urlcolor     = blue, 
	linkcolor    = blue, 
	citecolor   = red 
}

\newtheorem{problem}{Problem}[section]

\newtheorem{theorem}{Theorem}[section]
\newtheorem{lemma}{Lemma}[section]

\newtheorem{proposition}{Proposition}[section]
\newtheorem{corollary}{Corollary}[section]

\newtheorem{definition}{Definition}[section]

\def\sphere{S^{n-1}}
\def\R{\mathbb{R}}
\def\cK{\mathscr{K}}  
\def\deV{\widetilde{C}_{\Theta}}
\def\depsi{\Psi}
\def\ball{B^n}
\def\Rn{{\mathbb{R}^n}}
\def\cKn{\mathscr{K}_o^n}
\def\cKo{\mathscr{K}_{(o)}^n}
\def\N{\mathbb{N}}
\def\cH{\mathscr{H}}
\def\dvetV{\widetilde{V}_G}
\def\cG{\mathscr{G}}

\def\cC{\mathcal{C}}
\def\cM{\mathcal{M}}
\def\cKe {\mathscr{K}_e^n}

\def\cKen{\cK_e^n}
\def\polar{K^*}
\def\wp{G}
\def\VGL{\widetilde{V}_{G,\lambda}}
\def\cB{\mathcal{B}}
\def\cL{\mathcal{L}}

\def\bt{\begin{theorem}}
	\def\et{\end{theorem}}
\def\bl{\begin{lemma}}
	\def\el{\end{lemma}}
\def\bp{\begin{proposition}}
	\def\ep{\end{proposition}}
\def\bd{\begin{definition}}
	\def\ed{\end{definition}}
\def\bc{\begin{corollary}}
	\def\ec{\end{corollary}}

\newcommand{\conv}{{\mathrm{conv}}\,}

\begin{document}
	
	\title{On the Musielak-Orlicz-Gauss image problem}
	\author{  Qingzhong Huang, Sudan Xing, Deping Ye and Baocheng Zhu}
	
	\date{}
	\maketitle

	\begin{abstract}
		In the present paper we initiate the study of the Musielak-Orlicz-Brunn-Minkowski theory for convex bodies. In particular, we develop the Musielak-Orlicz-Gauss image problem aiming to characterize the Musielak-Orlicz-Gauss image measure of convex bodies. For a convex body $K$, its Musielak-Orlicz-Gauss image measure, denoted by $\widetilde{C}_{\Theta}(K, \cdot)$, involves a triple $\Theta=(G, \depsi, \lambda)$ where $G$ and $\Psi$ are two Musielak-Orlicz functions defined on $S^{n-1}\times (0, \infty)$ and $\lambda$ is a nonzero finite Lebesgue measure on the unit sphere $S^{n-1}$. Such a measure can be produced by a variational formula of  $\widetilde{V}_{G, \lambda}(K)$ (the general dual volume of $K$ with respect to $\lambda$) under the perturbations of $K$ by the Musielak-Orlicz addition defined via the function $\Psi$. The Musielak-Orlicz-Gauss image problem contains many intensively studied Minkowski type problems and the recent Gauss image problem as its special cases. Under the condition that $G$ is decreasing on its second variable, the existence of solutions to this problem is established.

	\end{abstract}

	\section{Introduction}\setcounter{equation}{0}
	The groundbreaking work \cite{LYZActa} by Huang, Lutwak, Yang and Zhang provides an extraordinarily beautiful connection between the Brunn-Minkowski theory for convex bodies and its dual. Among the most elegant concepts in \cite{LYZActa} are the  $q$-th dual curvature measures. These measures not only give the conceptual dual of the Federer's curvature measures, but also can be derived from the first order variation of the $q$th dual volume under the logarithmic perturbations of given convex bodies. Let $\cKo$ be the set of all convex compact subsets in $\Rn$ containing the origin in their interiors. By the $q$th dual volume of $K\in \cKo$, we mean  $$\widetilde{V}_q(K)=\frac{1}{n}\int_{\sphere}\rho_K^q(\xi)\,d\xi,$$ where $0\neq q\in \R$, $\,d\xi$ is the canonical spherical measure on $\sphere$ and $\rho_K: \sphere \rightarrow [0, \infty)$ is the radial function of $K$. By the logarithmic perturbations of $K\in \cKo$, we mean the family of convex bodies $[h_K\cdot e^{\epsilon g}]\in \cKo$, the  Wulff shapes generated by  $h_K\cdot e^{\epsilon g}$ (see  \eqref{july287} for the definition of the Wulff shape), where  $\epsilon\in \R$ is small enough, $h_K: \sphere \rightarrow [0, \infty)$ is the support function of $K$,  and $g: \sphere \rightarrow \R$ is a continuous function. Regarding the $q$-th dual curvature measures is the remarkable dual Minknowski problem \cite{LYZActa}: {\em given a real number $q$ and a nonzero finite Borel measure $\mu$ defined on  $\sphere$, can one find a  $K\in \cKo$ so that $\mu=\widetilde{C}_q(K, \cdot)$, with $\widetilde{C}_q(K, \cdot)$ the $q$-th dual curvature measure of $K\in \cKo$?} Since its introduction, the dual Minkowski problem has received a lot of attention, see e.g., \cite{BHP, BLYZZ2017, ChenLi-2018, Henk, Huangjiang, JiangWu, LiShengWang, WangFangZhou,  WangZ,  zhao, zhao-jdg}.

	The dual Minkowski problem has been pushed forward to the $L_p$ dual Minkowski problem \cite{LYZ-Lp} by Lutwak, Yang and Zhang and to the general dual Orlicz-Minkowski problem \cite{GHWXY,GHXY} by Gardner  {\it et al}. The latter one asks: {\em given  two continuous functions $G: \sphere \times (0, \infty)\to \R$ and  $\psi: (0, \infty)\to (0,\infty)$, under what conditions on a nonzero finite Borel measure $\mu$ defined on  $\sphere$ do there exist a $K\in \cKo$ and a constant $\tau\in \R$ so that $\mu=\tau \widetilde{C}_{G,\psi}(K, \cdot)$?} Here,  $\widetilde{C}_{G,\psi}(K, \cdot)$ is  the general dual Orlicz curvature measure for $K\in \cKo$ and can be formulated by: for any Borel set $\omega \subset \sphere$,   \begin{align} \label{c-g-21-1-10} \widetilde{C}_{G,\psi}(K, \omega)=\frac{1}{n}\int_{\pmb{\alpha}^*_K(\omega)} \frac{\rho_{K}(\xi)\, \wp_t(\xi, \rho_K(\xi)) }{\psi(h_{K}(\alpha_K(\xi)))}\,d\xi, \end{align} where $\alpha_K$ is the radial Gauss image of $K$, $\pmb{\alpha}^*_K$ is the reverse radial Gauss image of $K$ (see Section \ref{section-2} for detailed information), and $G_t$ is the first order partial derivative of $G$ with respect to its second variable.  Like the $q$-th dual curvature measure, the measure  $\widetilde{C}_{G,\psi}(K, \cdot)$  for $K\in \cKo$ can be obtained via the first order variation of the general dual volume $$\widetilde{V}_G(K)=\int_{\sphere}G(\xi, \rho_K(\xi))\,d\xi$$ in terms of the Orlicz $L_{\varphi}$ addition $\varphi^{-1}[\varphi(h_K)+\epsilon \varphi(g)],$ where $\varphi: (0, \infty)\to \R$ is a strictly monotonic function whose first order derivative $\varphi'$ satisfies $\psi(t)=t\varphi'(t)$ for $t\in (0,\infty)$.  

One of the biggest advantages for the general dual Orlicz-Minkowski problem \cite{GHWXY,GHXY} is its power to integrate various Minkowski type problems into a unified formula.    Here we give some special cases. First of all, if $G=t^n/n$,  then the $L_p$ and logarithmic Minkowski problems \cite{Lu93, min1897, min1903} are related to $\varphi(t)=t^p$ for $0\neq p\in \R$ and $\varphi(t)=\log t$, respectively. These problems have great impact on the development of the $L_p$ Brunn-Minkowski theory for convex bodies and have received immense attention, see   \cite{BHZ2016,BLYZ2013,  CLZ-2019, chen, chouw06, HLW-1, HugLYZ, JLW-1,  JLZ2016,  LuWang-1, Lut-Oli-12, LYZ04, ShengYi,  Umanskiy,  zhug20141, zhug2015b,GZhu2015II,  zhug2017} among others.   The Orlicz-Minkowski problem \cite{HLYZ2010} is related to the case when $G=t^n/n$ and $\varphi$ is a non-homogeneous function.  Solutions to the Orlicz-Minkowski problem  can be found in, e.g.,  \cite{BBColesanti, Bryanivakisch,   Huanghe2012, JianLu, liaijun2014,   LiuLu1, sun, SunLong,SunZhang,WuXiLeng2019, WuXiLeng2020}. When $G=t^q/n$ and $\varphi(t)=t^p$ for $0\neq p\in \R$, the general dual Orlicz-Minkowski problem  reduces to the $L_p$ dual Minkowski problem  \cite{LYZ-Lp}; contributions to this problem can be seen in, e.g., \cite{BorFod,  Chen-H-Z, ChenLi, ChenTWX1, HuangZhao, JWW-CVPDE-2021, LiLiuLu, ShengXia}. By letting $G(u, t)=\log t$ for all $(u, t)\in \sphere\times (0, \infty)$,  $\widetilde{V}_G(K)$ for $K\in \cKo$ reduces to the dual entropy of $K$; in this case one can get the ($L_p$ and Orlicz) Aleksandrov problems  \cite{Alexs1942, FH, HLYZ} (see also \cite{LiShengWang, zhao-pams}). Lastly,  the general dual Orlicz-Minkowski problem also extends the dual Orlicz-Minkowski problems \cite{XY2017-1, ZSY2017} and the Minkowski problem for Gaussian measures \cite{HuangXiZhao}.   Solutions to the (general) dual Orlicz-Minkowski problem by using the techniques from partial differential equations can be found in \cite{ChenKuyLuXiang,ChenTWX, LiuLu}.

In view of formula \eqref{c-g-21-1-10},  it is the radial Gauss image $\pmb{\alpha}_K: \sphere \to \sphere$ (or more precisely, the reverse radial Gauss image $\pmb{\alpha}_K^*$)  which plays an essential role to transfer $\,d\xi$ to   $\widetilde{C}_{G,\psi}(K, \cdot)$. Considering the importance of the reverse radial Gauss image, a new innovative problem bearing the flavour of the Minkowski type problems has been proposed in a recent paper  \cite{BLYZZ2020}  by   B\"{o}r\"{o}czky \emph{et al}. Such an elegant problem was named as the Gauss image problem which asks: {\it under what conditions on two given spherical Borel measures $\lambda$ and $\mu$, does there exist a $K\in \cKo$ such that $\mu=\lambda(\pmb{\alpha}_K(\cdot))$?}  (See Problem \ref{Gauss-I-p} for more general version.) The Gauss image problem involves two pre-given measures, and this is a major difference from the Minkowski type problems requiring only one pre-given measure.  As mentioned in \cite{BLYZZ2020}, if $\,d\lambda(\xi)=\,d\xi$, the Gauss image problem  reduces to a Minkowski type problem. We would like to comment that if $\, d\lambda(\xi)=p_{\lambda}(\xi)\,d\xi$ with a continuous function $p_{\lambda}: \sphere\to (0, \infty)$, the Gauss image problem indeed becomes a special case of the general dual Orlicz-Minkowski problem \cite{GHWXY,GHXY}. However, if the measure $\lambda$ does not have a continuous density with respect to $\,d\xi$ or even is not absolutely continuous with respect to $\,d\xi$, then the Gauss image problem is different from the Minkowski type problems. Under some mild conditions on $\lambda$ and $\mu$, the existence and uniqueness of solutions to the Gauss image problem have been established in \cite{BLYZZ2020}. See \cite{ChenWX} for smooth solutions to the Gauss image problem.

The present paper has two major goals. The first one is to provide a unified formulation to integrate the Minkowski type problems and the Gauss image problem.  The second one is to further push forward these problems to their next generation; this has the potential to initiate a brand-new theory: the  Musielak-Orlicz-Brunn-Minkowski theory for convex bodies. A closer observation on the Gauss image problem and the general dual Orlicz-Minkowski problem indicates that a triple $\Theta=(G, \Psi, \lambda)$ containing  three parameters shall be needed to fulfill these goals. Here $G$ and $\Psi$ are two Musielak-Orlicz functions defined on $\sphere\times (0, \infty)$ and $\lambda$ is a  spherical Lebesgue  measure on $\sphere$.  The function $G$ and the measure $\lambda$ are used to define the general dual volume of $K\in \cKo$ with respect to $\lambda$, namely, \begin{align*}
        \VGL(K)=\int_{\sphere} \wp(\xi, \rho_K(\xi))\, d\lambda(\xi).
        \end{align*}  The strictly monotone function $\depsi$ is used to define the addition of functions which produces the perturbation of convex bodies. In fact, the $L_{\depsi}$ addition of continuous functions $f$ and $g$ on $\Omega\subseteq \sphere$ can be formulated by 
        \begin{align}\label{genplus21-21-1-1} \depsi(\xi, f_{\varepsilon}(\xi))=\depsi(\xi, f(\xi))+\varepsilon g(\xi),\ \ \ \xi\in \Omega,
\end{align} where $\varepsilon\in \R$ is small enough and $f$ is strictly positive.  

In Section \ref{Section-var}, we will  define the Musielak-Orlicz-Gauss image measure and proves a variational formula to derive this measure.  Let $\Theta=(G, \Psi, \lambda)$ be a given triple with  $\lambda$ a nonzero finite Lebesgue measure on $\sphere$, $G\in \cC$ where $\cC$ is defined in \eqref{G-def-2021-01-16},  and $\depsi\in \cC_I\cup \cC_d$ where $\cC_I$ and  $\cC_d$ are defined in \eqref{G-I-2021-01-16-I} and \eqref{G-d-2021-01-16-d}, respectively. The Musielak-Orlicz-Gauss image measure $\deV(K, \cdot)$ for $K\in \cKo$ is defined as follows (see Definition \ref{m---111}): for each Borel set $\omega\subseteq \sphere$,    \begin{align*}
\deV(K, \omega)=\int_{\pmb{\alpha}^{*}_K(\omega)} \frac{\rho_K(\xi) G_t(\xi, \rho_K(\xi))}{h_{K}(\alpha_K(\xi))\depsi_t(\alpha_K(\xi), h_{K}(\alpha_K(\xi)))}\,d\lambda(\xi),
\end{align*} where $G_t$ and $\depsi_t$ are the first order partial derivatives of $G$ and $\depsi$ with respect to their second variables.  Under certain additional conditions on the measure $\lambda$ and the set $\Omega\subseteq\sphere$, the following variational formula can be established in Theorem \ref{ovev-cor}:  
\begin{align*}
\lim_{\varepsilon\rightarrow 0}\frac{\widetilde{V}_{G,\lambda}([f_{\varepsilon}])-\widetilde{V}_{G,\lambda}([f])}{\varepsilon}&=
\int_{\Omega} g(u)\, d\deV([f], u),
\end{align*} where $[f]$ denotes the Wulff shape of $f$, and $f_{\varepsilon}$ is given in \eqref{genplus21-21-1-1}.

In Section \ref{M:4},  we will propose our  Musielak-Orlicz-Gauss image problem (i.e.,  Problem \ref{MOGIP}): {\em  Under what conditions on $\Theta=(G, \Psi, \lambda)$ and a nonzero finite Borel measure $\mu$ on $\sphere$ do there exist a $K\in \cKo$ and a constant $\tau\in \R$ such that  $\mu =\tau \deV(K,\cdot)?$}  As one can see in Section \ref{M:4}, the Musielak-Orlicz-Gauss image problem extends all previously mentioned Minkowski type problems and the Gauss image problem in their (arguably) most general formulations. Some special cases of particular interest are discussed; these include the Musielak-Orlicz-Minkowski problem  (i.e., Problem \ref{MOMP-vol}) when $ \VGL(\cdot)$ is the volume, the  dual Musielak-Orlicz-Minkowski problem (i.e., Problem \ref{MOMP-dual}) when $\,d\lambda(\xi)=\,d\xi$, and the Musielak-Orlicz-Aleksandrov problem (i.e., Problem \ref{MOMP-Alex}) when $G=\log t$ and $\,d\lambda(\xi)=\,d\xi$.  As byproducts, we obtain some fully nonlinear Monge-Amp\`{e}re partial differential equations. Indeed,  if $\mu$ and $\lambda$ have continuous density functions $p_{\mu}$ and $p_{\lambda}$, respectively,  then  the Musielak-Orlicz-Gauss image problem could be reformulated by \begin{align*}
p_{\mu}= \tau  \frac{P(\bar{\nabla}h+h\iota) \,\det(\bar{\nabla}^2h+hI)}{\depsi_{t}(\cdot, h)}  p_{\lambda}\! \left(\!\frac{\bar{\nabla} h+h \iota}{|\bar{\nabla} h+h \iota|}\!\right), \end{align*}
where  $\tau\in \R$, $P(y)=|y|^{1-n}G_{t}(\frac{y}{|y|}, |y|)$ for $y\in \R^n$ with $|y|$ the Euclidean norm of $y\in \R^n$,  $\iota$ denotes the identity map on $\sphere$,      $\bar{\nabla}h$ and $\bar{\nabla}^2h$  are the gradient and the Hessian matrix of $h$ with respect to an orthonormal frame on $\sphere$, and $I$ is the identity matrix. 

 Under the condition that $G$ is strictly decreasing on its second variable,  the existence of solutions to the Musielak-Orlicz-Gauss image problem will be established in Sections \ref{solution-8-25} and \ref{even-8-25}. (Section \ref{even-8-25} deals with the Musielak-Orlicz-Gauss image problem for even data). A typical result is Theorem \ref{solution-general-dual-Orlicz-main theorem-11-270}, which is stated below.  

\vskip 2mm \noindent {\bf Theorem \ref{solution-general-dual-Orlicz-main theorem-11-270}.} {\em  Let $\lambda$ and $\mu$ be two nonzero finite Borel measures on $S^{n-1}$ that are not concentrated on any closed hemisphere. Assume that  $\lambda$ is absolutely continuous with respect to  $\,d\xi$. Then, there exists a $K\in \cKo$  such that  
    \begin{align}\label{msol-8-5-1-intro}
    \frac{\mu}{|\mu|}=\frac{\deV(K,\cdot)}{\deV(K, \sphere)},
    \end{align} if $G$ and $\depsi$ satisfy one of the following conditions: 
   
   \vskip 2mm \noindent    i) $\wp\in \cG_d$ with  $\cG_d$ given by \eqref{G-2021-1-16}, and $\depsi\in \cC_I$ with $\cC_I$ given by \eqref{G-I-2021-01-16-I} such that 
 $$\lim_{t\rightarrow \infty} \depsi(\xi, t)=+\infty \ \ \mathrm{for\ each}\  \xi\in \sphere;$$ 
 
   \vskip 2mm  \noindent  ii) $\depsi\in \cG_I$ with $\cG_I$ given by \eqref{G-2021-1-16}, and $G\in \cC_d$ with $\cC_d$ given by \eqref{G-d-2021-01-16-d} such that $$\lim_{t\rightarrow 0^+} G(\xi, t)=+\infty \ \ \mathrm{for\ each}\ \xi\in \sphere. $$     
   If, in addition, $\lambda$ is strictly positive on nonempty open subsets of $\sphere$, then  the assumption on $\mu$, i.e., $\mu$ is a nonzero finite Borel measure on $S^{n-1}$ that is not concentrated on any closed hemisphere, is also necessary for  \eqref{msol-8-5-1-intro} holding true for some $K\in \cKo$.}  
   
 Under special choices of $\Theta=(G, \depsi, \lambda)$, one also obtains the solutions to  the  dual Musielak-Orlicz-Minkowski problem  and the Musielak-Orlicz-Aleksandrov problem. 
 We would like to mention that, under additional condition on $\mu$ (i.e.,  $\mu$ vanishes on great subspheres),  the Musielak-Orlicz-Gauss image problem for even data can also be solved when $\depsi\in \cC_d$ and  $\wp\in \cG_d$ (see Theorem \ref{rl-even-ee}).  Corollary \ref{MMM-2021-1-16} provides the existence of solutions to the Musielak-Orlicz-Aleksandrov problem for even data, under both $\depsi\in  \cG_I$  and $\depsi \in \cG_d$.   Solutions to the Musielak-Orlicz-Gauss image problem related to ``an increasing function" $G\in \cC_I$ will be studied in our future work \cite{HXYZ-2}; while the solutions to the Musielak-Orlicz-Gauss image problem by the technique of flows will be provided in \cite{LSYY-2021}.

\section{Preliminaries and notations}\label{section-2}\setcounter{equation}{0}

Denote by $\N$ the set of all natural numbers. Let $n\in \N$ be such that $n\geq 2$. In the  $n$-dimensional Euclidean space $\Rn$, let $\sphere$ be the  unit sphere  and $\ball$ be the unit Euclidean ball of $\Rn$; namely, $$\sphere=\{x\in\Rn: |x|=1\} \ \ \mathrm{and} \ \ \ball= \{x\in
\Rn:|x|\leq 1\},$$  where $|x|$ denotes the Euclidean norm of $x\in \R^n$.  The origin of $\Rn$ is denoted by $o$, and the inner product of $x, y
\in \Rn$ is written by  $x\cdot y$.  For $x\neq o$, let $\overline{x}=x/|x|.$ Denote by $\cH^{n-1}$ the $(n-1)$-dimensional Hausdorff measure. In particular, let $\,d\xi$ be the canonical spherical Lebesgue measure on $\sphere$.

Let $K\subseteq \Rn$ be a nonempty, compact and convex set. Define $h_K: \R^n\rightarrow \R$, the support function of $K$, to be \begin{align}\label{support-function--111}
h_K(x)=\max \big\{x\cdot y: y\in K\big\} \ \ \ \mathrm{for} \ x\in \Rn. \end{align}  It is easily checked that $h_K(ru)=rh_K(u)$ for $r>0$ and $u\in\sphere$. For two  nonempty, compact and convex sets $K, L\subseteq \Rn$, let $$d_H(K, L)=\max_{u\in \sphere}|h_K(u)-h_{L}(u)|.$$  The convergence of $K_i\rightarrow K$ in the Hausdorff metric, where $K_i, K\subseteq \Rn$ for all $i\in \N$ are nonempty, compact and convex sets, is defined by $$\lim_{i\rightarrow \infty}d_H(K_i, K)=0. $$ 
 
A convex body is a compact and convex subset of $\R^n$ whose interior is nonempty.  The  Blaschke selection theorem asserts that every sequence of convex bodies, if uniformly bounded, must have a subsequence converging to a compact convex set in $\Rn$.   Let $\cKn$ be the set of all convex bodies containing the origin. For $K\in \cKn$, let   $\rho_K:
\mathbb{R}^n\setminus\{o\}\rightarrow[0,\infty)$ stand for the radial function of $K$. That is,  
\begin{align}\label{support-function--222}
   \rho_K(x)=\max \big\{\lambda\ge
0:\lambda x\in K\big\},
\end{align}  for   $x\in \mathbb R^n\backslash\{o\}$.  It is easily checked that  $\rho_K(ru)=r^{-1}\rho_K(u)$ for $r>0$ and $u\in\sphere$.    In this paper, we are mainly interested in the class of convex bodies $\cKo$ which consists of all convex bodies whose interiors contain the  origin $o$.  We say $K\in \cKo$ is origin-symmetric if $-K=K$, where $aK=\{ax:  x\in K\}$ for $a\in \R$. The subclass $\cKen$ of   $\cKo$ denotes the set of all origin-symmetric convex bodies.   For $K\in \cKo$, let $$
    \polar=\big\{x\in\mathbb R^n: x\cdot y\le 1\quad\textrm{for all}~y\in
    K\big\}.
$$ It can be easily checked that $\polar \in \cKo$ for $K\in \cKo$ having the following properties: \begin{align}\label{bi-polar--12}
    h_{\polar}(u) \rho_K(u) =1 \quad \textrm{and}\quad
    \rho_{\polar}(u) h_K(u)=1,
\end{align} for $u\in \sphere$. The convex body $\polar$ is called the polar body of $K\in \cKo$.

Formulas \eqref{support-function--111} and \eqref{support-function--222} naturally bring many key concepts which play important roles in this paper. The first one is the so-called radial map  $r_K$ of $K\in \cKo$ which maps $u\in \sphere$ to $r_K(u)=\rho_K(u)u\in \partial K$, the boundary of $K$. The map $r_K$ is invertible for $K\in \cKo$ and its reverse, denoted by $r_K^{-1}$, maps $x\in \partial K$ to $\sphere$ by letting $r_K^{-1}(x)=\overline{x}.$ The second one is the Gauss map $\pmb{\nu}_K: \partial K\mapsto \sphere$ which maps an $x\in \partial K$ to all $u\in \sphere$ satisfying $x\cdot u =h_K(u)$.  Note that  $\pmb{\nu}_K$ may not be injective on $\partial K$. Let $$\sigma_K=\big\{x\in \partial K: \  \pmb{\nu}_K(x) \ \mathrm{contains\ two\ or\ more \ elements} \big\}\subseteq \partial K.$$ It is well known that $\cH^{n-1}(\sigma_K)=0$ \cite[p.84]{Sch}. Clearly, $\pmb{\nu}_K$ is injective in $\mathrm{reg}(K)=\partial K\setminus \sigma_K$; in this case, for simplicity, we write $\nu_K(x)$ for $\pmb{\nu}_K(x)$ if  $x\in \mathrm{reg}(K)$.   
Likewise, the inverse Gauss map $\pmb{\nu}^{-1}_K: \sphere \to \partial K$ maps $u\in \sphere$ to all $x\in \partial K$ such that $ x\cdot u =h_K(u).$ Note that $\pmb{\nu}^{-1}_K$ may not be injective and it is injective in the set $\sphere \setminus \eta_K$ where $$\eta_K=\big\{u\in \sphere:  \  \pmb{\nu}^{-1}_K(u) \ \mathrm{contains\ two \ or \ more\ elements} \big\}\subseteq \sphere.$$ 
Again $\cH^{n-1}(\eta_K)=0$ as shown in \cite[Theorem 2.2.11]{Sch}. Let $\nu^{-1}_K(u)=\pmb{\nu}^{-1}_K(u)$  for $u\in \sphere \setminus \eta_K$.  

Let $K\in \cKo$. One can define the radial Gauss image $\pmb{\alpha}_K: \sphere\rightarrow \sphere$ which maps $u\in \sphere$ to the set $\pmb{\alpha}_K(u)=\pmb{\nu}_K(r_K(u)).$ Namely,
$\pmb{\alpha}_K(u)$ is the set of all outer unit normal vectors of $\partial K$ at the point $\rho_K(u)u\in \partial K$. Define 
$$\omega_K=\big\{u\in \sphere \ \mathrm{such\ that}\  \pmb{\alpha}_K(u)\ \mathrm{contains\ two \ or \ more\ elements}  \big\}.$$ Note that $\cH^{n-1}(\omega_K)=0$ as shown in \cite[Theorem 2.2.5]{Sch} (or see  \cite[p.340]{LYZActa}). Let $ \alpha_K(u)=\pmb{\alpha}_K(u)$ for $u\in \sphere\setminus \omega_K$. 
 On the other hand, one can define the reverse radial Gauss image $\pmb{\alpha}^*_K: \sphere \to \sphere$ which maps $u\in \sphere$ to the set  $\pmb{\alpha}^*_K(u)=r_K^{-1}(\pmb{\nu}^{-1}_K(u)).$ Moreover, $\pmb{\alpha}^*_K$ is injective on the set $\sphere\setminus \eta_K$, and in this case, $\pmb{\alpha}^*_K(u)$ will often be written as $\alpha^*_K(u)$. According to \cite[Lemma 2.5]{LYZActa}, the radial Gauss image and its reverse can be connected through the polar body:
 for any $K\in \cKo$ and $\eta\subseteq \sphere$, one has \begin{align}\label{star}
\pmb{\alpha}^{*}_K(\eta)=\pmb{\alpha}_{K^{*}}(\eta).
\end{align}

Let $\cB$ and $\cL$ be the $\sigma$-algebras of spherical Borel and Lebesgue measurable subsets of $\sphere$,  respectively. 
We say $\lambda$ a spherical Lebesgue submeasure if $\lambda: \cL \to [0, \infty)$ satisfies that $\lambda(\omega)=0$ if $\omega$ is an empty set;  $\lambda(\omega_1)\leq \lambda(\omega_2)$ for all $\omega_1, \omega_2\in \cL$ such that $\omega_1\subseteq \omega_2$; and 
 $\lambda(\bigcup_{i=1}^{\infty} \omega_i)\leq \sum_{i=1}^{\infty}  \lambda(\omega_i)$ if $\omega_i\in \cL$ for all $i\in \N$.  A spherical Borel submeasure can be defined in a similar way with $\cL$ replaced by $\cB$.   A nonzero finite Borel measure $\mu$ on $\sphere$ is said to be not concentrated on any closed hemisphere of $\sphere$, if for any $u\in \sphere$, one has \begin{align}\label{not-concentration-1} \int_{\sphere} (u\cdot \xi)_+\,d\mu(\xi)>0\end{align} with $a_+=\max\{a, 0\}$ for $a\in \R$. 
This is equivalent to  the following statement: $\mu(\sphere\setminus\omega)>0$ if $\omega$ is  an arbitrary closed hemisphere of $\sphere$. 
 
  Let $K\in \cKo$.  The surface area measure $S(K, \cdot)$ is defined by \begin{align}\label{surface-8-7} S(K, \omega)=\cH^{n-1}(\pmb{\nu}_K^{-1}(\omega)
  )\ \ \ \mathrm{for} \ \omega\in \cB.\end{align}  Similarly, the composition of a spherical Lebesgue submeasure $\lambda$ and the reverse radial Gauss image $\pmb{\alpha}^*_K$ naturally defines a spherical Borel submeasure on $\sphere$ (see \cite{BLYZZ2020}); such a spherical Borel submeasure is named as the reverse Gauss image measure of $\lambda$ via $K\in \cKo$ and is denoted by $\lambda^*(K, \cdot)$. That is,  for each $\omega\in \cB$, \begin{align}\label{r-G-i-measure} 
  \lambda^*(K, \omega)=\lambda(\pmb{\alpha}_K^*(\omega))=\lambda(\pmb{\alpha}_{\polar}(\omega)). 
  \end{align} The Gauss image measure of $\lambda$ via $K$ \cite{BLYZZ2020},  another spherical Borel submeasure, can be defined by \begin{align}\label{G-i-measure} 
  \lambda(K, \omega)=\lambda(\pmb{\alpha}_K(\omega)). 
  \end{align} According to \eqref{r-G-i-measure} and \eqref{G-i-measure}, it is easily checked that  $ \lambda^*(K, \cdot)= \lambda(\polar, \cdot)$  for each $K\in \cKo$.  It has been proved in \cite[Lemma 3.3]{BLYZZ2020} that if $\lambda$ is a Borel measure which is absolutely continuous with respect to $\,d\xi$, then $\lambda^*(K, \cdot)$ for $K\in \cKo$ is a spherical Borel measure; moreover,  for bounded Borel function $f$ defined on $\sphere$, one has  \begin{align} \int_{\sphere} f(u)\,d\lambda^* (K, u)&=\int_{\sphere} f(\alpha_K(\xi))\,d\lambda(\xi) \nonumber;   \\   \int_{\sphere} f(u)\,d\lambda (K, u)&=\int_{\sphere} f(\alpha^*_{K}(\xi))\,d\lambda(\xi). \label{int-G-i-measure-1-1} \end{align}
   
   When $\,d\lambda(\xi)=\,d\xi$, $\lambda(K, \cdot)$ for $K\in \cKo$ reduces to the Aleksandrov's integral curvature \cite{Alexs1942}, which will be denoted by $J(K, \cdot)$.  Similarly, one can let $J^*(K, \cdot)=J(\polar, \cdot)$ for $K\in \cKo$.  Hence,   the following formulas hold (see \cite{HLYZ}):   
   \begin{align*} \int_{\sphere} f(u)\,dJ^* (K, u)=\int_{\sphere} f(\alpha_K(\xi))\,d\xi \ \ \mathrm{and}\  \ \int_{\sphere} f(u)\,dJ (K, u)=\int_{\sphere} f(\alpha^*_{K}(\xi))\,d\xi.  \end{align*}

  We shall need the following lemma, which is essentially a restatement of  \cite[Lemmas 5.3 and
5.4]{BLYZZ2020}. We will present a proof
here for completeness.

 \bl\label{conc} Let
\(\lambda\) be an absolutely continuous Borel measure on \(S^{n-1}\)
that is strictly positive on nonempty open subsets of $\sphere$. Then, 
$\lambda(K, \cdot)$  for $K\in
\cKo$  is not
concentrated on any closed hemisphere. In particular, for any $u \in \sphere$, one has \begin{align}\label{not-concentration-1-22} \int_{\sphere} (u\cdot v)_+\,d\lambda(K, v)=\int_{\sphere}(u\cdot \alpha_K^*(\xi))_+\,d\lambda(\xi)>0.\end{align} 
\el
        \begin{proof}    Let $\omega\subseteq \sphere$ be a nonempty subset. Define \begin{align*} \mathrm{cone} (\omega)=\big\{ tu: \ u\in \omega \ \mathrm{and}\  t\geq 0\big\}\ \ \mathrm{and}\ \  \omega^*=\big\{v\in \sphere: u\cdot v\leq 0 \ \mathrm{for} \ u\in \omega\big\}.\end{align*}   If $\mathrm{cone} (\omega)$ is a proper convex subset of $\Rn$, then $\omega$ is called a spherically convex set; in this case,  $\omega$ certainly is contained in a closed hemisphere of $\sphere$. 
        
           It is easily checked that $\pmb{\alpha}_K(\omega)\subseteq \sphere\setminus \omega^*$ holds for all spherically convex subset $\omega\subseteq \sphere$; this can be seen from \cite[Lemma 3.2]{BLYZZ2020} where the following statement is also given: $(\sphere\setminus \omega^*)\setminus \pmb{\alpha}_K(\omega)$ has interior points. It follows from  \eqref{G-i-measure} and the assumption on $\lambda$ (in particular, the positiveness of $\lambda$ on nonempty open subsets of $\sphere$) that \begin{align}\label{2021-02-01-1} \lambda(K, \omega)= \lambda\left(\pmb{\alpha}_{K}(\omega)\right)<\lambda\left(\sphere \backslash \omega^*\right)
        \end{align} holds for any spherically convex $\omega\subseteq \sphere$; this argument can be seen in  \cite[Lemma 3.7]{BLYZZ2020}. 
     
     Note that   $\lambda(K, \sphere)=\lambda\left(\boldsymbol{\alpha}_{K}(\sphere)\right)=\lambda\left(\sphere\right)$. Then, for any spherically convex set $\omega\subseteq \sphere$, by  \eqref{2021-02-01-1}, one gets 
    \begin{align*}\lambda\left(\omega^*\right) =\lambda\left(\sphere\right) - \lambda\left(\sphere \setminus 
    \omega^*\right)  <\lambda\left(\sphere\right)  -\lambda(K, \omega)=\lambda(K, \sphere\setminus\omega).
        \end{align*} Now let $\omega$ be an arbitrary closed hemisphere of $\sphere$. Then $\omega^*$ contains only one vector in $\sphere$ and  $\sphere\setminus \omega$ is an open hemisphere of $\sphere$. Consequently, $\lambda(K, \cdot)$ is not concentrated on any closed hemisphere because $\lambda(K, \sphere\setminus\omega)>\lambda\left(\omega^*\right) \geq 0. $ 
        Formula \eqref{not-concentration-1-22}  is then an immediate consequence of \eqref{not-concentration-1} and \eqref{int-G-i-measure-1-1}.  \end{proof}

Regarding the Gauss image measure, the following Gauss image problem has been posed in \cite{BLYZZ2020}. 
\begin{problem}[The Gauss image problem]\label{Gauss-I-p}  Let $\lambda$ be a spherical Lebesgue submeasure and $\mu$ be a spherical Borel submeasure. Under what conditions on $\lambda$ and $\mu$, does there exist a $K\in \cKo$ such that $\mu(\omega)=\lambda(K,\omega)$ holds for all $\omega\in \cB$?  
\end{problem}  

Solutions to  the Gauss image problem can be found in  \cite{BLYZZ2020, ChenWX}.  When $\,d\lambda(\xi)=\,d\xi$, Problem \ref{Gauss-I-p} becomes the  classical Aleksandrov problem aiming to characterize the Aleksandrov's integral curvature.  As mentioned in the introduction, our goal in this paper is to introduce a problem extending the Minkowski type problems and the  Gauss image problem in their (arguably) most general setting by taking use of the Musielak-Orlicz functions. The precise definition for the Musielak-Orlicz functions can be found in  e.g., \cite{Harjulehto, Musielak},  and they play important roles in the analysis of the  Musielak-Orlicz space (or the generalized Orlicz space). 

The  Musielak-Orlicz  functions in e.g., \cite{Harjulehto, Musielak}, are usually referred to strictly positive and nondecreasing functions.  However, throughout this paper, when we refer to the Musielak-Orlicz  functions, we mean  $G\in \cC$ with \begin{align} \cC=\big\{G: \sphere\times(0, \infty)\to \R  \ \mathrm{such\ that}\ G  \ \mathrm{and}\  G_t\ \mathrm{ are\ continuous\ on}\ \sphere\times(0, \infty)\big\},   \label{G-def-2021-01-16}   \end{align}  where   $G_t$ denotes the partial derivative of $G$ with respect to the second variable, namely, $$G_t(\xi, t)=\frac{\partial G}{\partial t}(\xi, t) \ \ \mathrm{for}\ \ (\xi, t)\in \sphere\times (0, \infty).$$ As the function class $\cC$ does contain all (smooth enough) Musielak-Orlicz  functions
defined in, e.g., \cite{Harjulehto, Musielak}, taking use of the function class $\cC$ not only provides convenience in later context but also gives reasons to name $\deV(K, \cdot)$ the Musielak-Orlicz-Gauss image measure (see Definition \ref{m---111}). 

Let $\cG_I$ and $\cG_d$ be subclasses of $\cC$ defined by   \begin{align}
\cG_I = \big\{G\in \cC:  G  \ \mathrm{satisfies \ condition\ ({\bf A})} \big\}  \  \mathrm{and}\   \cG_d = \big\{G\in \cC:  G  \ \mathrm{satisfies \ conditions\ ({\bf B})} \big\},  \label{G-2021-1-16}
\end{align} where conditions ({\bf A}) and ({\bf B}) are given below:  \begin{itemize}  \item[({\bf A})]   $G: \sphere\times (0, \infty)\to (0, \infty)$ satisfies that, for each  $u\in \sphere$,   $G_t(u, \cdot)$  is strictly positive on $(0, \infty)$, $\lim_{t\to 0^+}G(u, t)=0$, and  $\lim_{t\to \infty}G(u, t)=\infty;$ 
 \item[({\bf B})] $G: \sphere\times (0, \infty)\to (0, \infty)$  satisfies that, for each $u\in \sphere$,  $G_t(u, \cdot)$  is strictly negative on $(0, \infty)$,  $\lim_{t\to 0^+}G(u, t)=\infty,$ and  $\lim_{t\to \infty}G(u, t)=0.$
 \end{itemize}  
The fact that $G\in \cG_I\cup \cG_d$ is assumed to be strictly positive is mainly for technique reasons and for convenience. Our arguments in later context mainly rely on the monotonicity of $G$ on its second variable, $\sup_{t>0}G(u, t)=+\infty$ for each $u\in \sphere$, and the fact that $G$ has controllable lower bounds;  our results in later context should still hold if the function $\inf_{t>0}G(u, t): \sphere \rightarrow \R$ is continuous on $\sphere$. 
 
We shall also need the following classes of functions: \begin{align} 
\cC_I&=\big\{G\in \cC:  \   G_t \ \mathrm{is\ strictly\ positive\ on\ } \sphere\times (0, \infty) \big\}, \label{G-I-2021-01-16-I}\\
\cC_d&=\big\{G\in \cC:  \  G_t \ \mathrm{is\ strictly\ negative \ on\ } \sphere\times (0, \infty)  \big\} \label{G-d-2021-01-16-d}.\end{align} When we write $G=\varphi(t)$ for some function $\varphi: (0, \infty)\rightarrow \R$, we mean $G(\xi, t)=\varphi(t)$ holding true for all $(\xi, t)\in \sphere\times (0, \infty).$  For $\depsi\in\cC$, let $\psi_{\xi}$ and  $\widetilde{\depsi}$ be the functions given by
 \begin{align}\label{formula220}
 \psi_{\xi}(t)=\depsi(\xi, t) \ \ \mathrm{and}\ \ \widetilde{\depsi}(\xi, t)=\depsi\big(\xi, 1/t\big)\ \ \mathrm{for}\ \ (\xi, t)\in \sphere\times (0, \infty).\end{align} 
Clearly,  $\widetilde{\depsi}_t(\xi, t)=-\frac{1}{t^2}\depsi_t(\xi, \frac{1}{t})$ and $\widetilde{\depsi}(\xi, t)\in \cC$. 
 
 \section{The  Musielak-Orlicz-Gauss image measure and related variational formula}\label{Section-var} \setcounter{equation}{0}
In this section,  the Musielak-Orlicz-Gauss image measure is introduced and a variational formula to derive such a measure is established. Hereafter,  $\lambda$ is always assumed to be a nonzero finite Lebesgue measure on $\sphere$.  For convenience, let    \begin{align*}  \cM&=\big\{\mathrm{nonzero\ finite\ Borel\ measures\ on\ }  \sphere \ \mathrm{ that\ are \ absolutely\ continuous\ w.r.t.}\  \,d\xi  \big\}. \end{align*}

 The definition for the Musielak-Orlicz-Gauss image measure is given below. 
\bd\label{m---111} Let $\Theta=(G, \Psi, \lambda)$ be a given triple with $G\in \cC$,  $\depsi\in \cC_I\cup \cC_d$, and $\lambda$ a nonzero finite Lebesgue measure on $\sphere$. Define   $\deV(K, \cdot)$,   the Musielak-Orlicz-Gauss image measure of $K\in \cKo$,  as follows: for each Borel set $\omega\in \cB$,                  \begin{align}\label{gencdef-7-1-115}
            \deV(K, \omega)=\int_{\pmb{\alpha}^{*}_K(\omega)} \frac{\rho_K(\xi) G_t(\xi, \rho_K(\xi))}{h_{K}(\alpha_K(\xi))\depsi_t(\alpha_K(\xi), h_{K}(\alpha_K(\xi)))}\,d\lambda(\xi). 
            \end{align}
        \ed

Let $\Theta_0=(\log t, \log t, \lambda)$. It follows from \eqref{r-G-i-measure} that,  for all $\omega\in \cB$, one has       \begin{align}\label{gencdef-7-1-115-uu}
           \widetilde{C}_{{\Theta}_0}(K, \omega)=\int_{\pmb{\alpha}^{*}_K(\omega)}  \,d\lambda(\xi)=\lambda(\pmb{\alpha}^{*}_K(\omega))=\lambda^*(K, \omega). \end{align} Formula \eqref{gencdef-7-1-115} implies that $ \deV(K, \cdot)$ is absolutely continuous with respect to $\lambda^*(K, \cdot)$, and  
      \begin{align}  \frac{\,d \deV(K, u)} {\,d \widetilde{C}_{{\Theta}_0}(K, u)} =\frac{\,d \deV(K, u)} {\,d\lambda^*(K, u)} = \frac{\rho_K(\alpha^*_K(u)) G_t(\alpha^*_K(u), \rho_K(\alpha^*_K(u)))}{h_{K}(u)\depsi_t(u, h_{K}(u))}\ \ \ \mathrm{for}\ \ u\in \sphere. \label{relation-G} \end{align}
        
When $\,d\lambda=\,d\xi$, \eqref{gencdef-7-1-115} reduces to 
       \begin{align}\label{gencdef-7-1-115--1}
            \widetilde{C}_{G, \depsi}(K, \omega)=\int_{\pmb{\alpha}^{*}_K(\omega)} \frac{\rho_K(\xi) G_t(\xi, \rho_K(\xi))}{h_{K}(\alpha_K(\xi))\depsi_t(\alpha_K(\xi), h_{K}(\alpha_K(\xi)))}\,d\xi. 
            \end{align}    In this case, if $\depsi=\varphi(t)$ for some function $\varphi: (0, \infty)\to \R$ whose derivative, denoted by $\varphi'$,  satisfies $t\varphi'(t)=\psi(t)$, then $\frac{1}{n}  \widetilde{C}_{G, \depsi}(K, \cdot)$ becomes the  general dual Orlicz curvature measure $\widetilde{C}_{G, \psi}(K, \cdot)$ in \cite[Definition 3.1]{GHWXY}.  Hence $\deV(K, \cdot)$ naturally extends $\widetilde{C}_{G, \psi}(K, \cdot)$ to its (arguably)  most  general setting;  and certainly contains  many well-known measures appeared in the Minkowski type problems as its special cases, including but not limited to the surface area measure \eqref{surface-8-7}, the $L_p$ surface area measure \cite{Lu93}, the Orlicz surface area measure  \cite{HLYZ2010}, the $L_{p}$  dual curvature measure \cite{LYZActa, LYZ-Lp}, the  dual Orlicz curvature measure \cite{XY2017-1, ZSY2017}, and  the Aleksandrov's integral curvature and its extensions \cite{Alexs1942, FH, HLYZ} (up to a difference of polarity of convex bodies).

 In the following definition, we define another measure which is closely related to  the Musielak-Orlicz-Gauss image measure.  It gives a great convenience to establish solutions to the Musielak-Orlicz-Gauss image problem.
\bd Let $\Theta=(G, \Psi, \lambda)$ be a given triple with $G\in \cC$,  $\depsi\in \cC_I\cup \cC_d$, and $\lambda$ a nonzero finite Lebesgue measure on $\sphere$.  Define $C_{\Theta}(K, \cdot)$,
the   polar  Musielak-Orlicz-Gauss image measure of $K\in \cKo$, 
  as follows: for each Borel set $\omega\in \cB$,  
\begin{align}\label{pd}
            C_{\Theta}(K, \omega)=\int_{\pmb{\alpha}^{*}_K(\omega)} \frac{\rho_K(\xi) G_t( \xi, \rho_K(\xi)) }{\rho_{K^*}(\alpha_K(\xi))\depsi_t
            (\alpha_K(\xi), \rho_{K^*}(\alpha_K(\xi)))}\,d\lambda(\xi). 
            \end{align} 
        \ed

Associated to $\Theta=(G, \depsi, \lambda)$,  let  $\widetilde{\Theta}=(G, \widetilde{\depsi},\lambda)$ with $\widetilde{\depsi}$ defined in (\ref{formula220}).  For $K\in \cKo$, one has  \begin{align}\label{wc}
            C_{\widetilde{\Theta}}(K, \cdot )=- \widetilde{C}_{\Theta}(K, \cdot). 
            \end{align} To this end,   by \eqref{bi-polar--12},
\eqref{gencdef-7-1-115}, \eqref{pd} and $\widetilde{\depsi}_t(\xi, t)=-\frac{1}{t^2}\depsi_t(\xi, \frac{1}{t})$, one gets, for any $\omega\in \cB$, 
\begin{align*}
 C_{\widetilde{\Theta}}(K, \omega)&=\int_{\pmb{\alpha}^{*}_K(\omega)} \frac{\rho_K(\xi) G_t(\xi, \rho_K(\xi)) }{\rho_{K^*}(\alpha_K(\xi))\widetilde{\depsi}_t (\alpha_K(\xi), \rho_{K^*}(\alpha_K(\xi)))}\,d\lambda(\xi)
\nonumber\\
&=-\int_{\pmb{\alpha}^{*}_K(\omega)} \frac{\rho_{K^*}(\alpha_K(\xi))\rho_K(\xi) G_t(\xi, \rho_K(\xi)) }{\depsi_t (\alpha_K(\xi), 1/\rho_{K^*}(\alpha_K(\xi)))}\,d\lambda(\xi)\nonumber\\
&=-\int_{\pmb{\alpha}^{*}_K(\omega)} \frac{\rho_K(\xi)
G_t(\xi, \rho_K(\xi))
}{h_{K}(\alpha_K(\xi))\depsi_t(\alpha_K(\xi), h_{K}(\alpha_K(\xi)))}\,d\lambda(\xi)
           \nonumber\\
&=- \widetilde{C}_{\Theta}(K, \omega).
\end{align*}

It is not hard to prove that  both  $\deV(K, \cdot)$ and $C_{\Theta}(K, \cdot)$,  for $K\in \cKo$, are finite signed Borel measures on $\sphere$.  The proof of this argument for $\deV(K,
\cdot)$ (and hence for $C_{\Theta}(K, \cdot)$ due to \eqref{wc})  is rather standard and follows from steps very similar to those in  \cite[p.9]{GHWXY} or \cite[p.351-352]{LYZActa}.  It can be also proved by using \eqref{relation-G}  and the fact that   $\lambda^*(K, \cdot)$ is a Borel measure on $\sphere$ (see \cite[Lemma 3.3]{BLYZZ2020}), thus the proof is omitted.       A standard argument,  based on the simple functions and a limit approach, shows that,  for any bounded Borel function $g:\sphere\to \R$,    
\begin{align} \int_{\sphere} g(u)\, d\deV(K, u) &= \int_{\sphere}\frac{g(\alpha_K(\xi))\rho_K(\xi) G_t(\xi, \rho_K(\xi))}
            {h_{K}(\alpha_K(\xi))\depsi_t(\alpha_K(\xi), h_{K}(\alpha_K(\xi)))}\,d\lambda(\xi),\label{new measue-11-27}\\ 
             \int_{\sphere} g(u)\, dC_{\Theta}(K, u)&=\int_{\sphere}\frac{g(\alpha_K(\xi))\rho_K(\xi) G_t(\xi, \rho_K(\xi)) }{\rho_{K^*}(\alpha_K(\xi))\depsi_t
            (\alpha_K(\xi), \rho_{K^*}(\alpha_K(\xi)))}\,d\lambda(\xi).\label{nn}
            \end{align}
 It is well known that, by letting $\xi=\overline{x}$ with $x\in \partial K$ for $K\in \cKo$ and  $\,dx=\,d\cH^{n-1}(x)$,  \begin{align}\label{vol-change-8-8} \int_{\sphere}f(\xi)\,d\xi=\int_{\partial K} (x \cdot \nu_K(x)) f(\overline{x})|x|^{-n}\,dx,\end{align}    (see \cite[(2.31)]{LYZActa}). Hence, if $\,d\lambda(\xi)=p_{\lambda}(\xi)\,d\xi$ with $p_{\lambda}: \sphere\rightarrow [0, \infty)$, then         
        \begin{align} \int_{\sphere} g(u)\, d\deV(K, u)  &= \int_{\partial K}g(\nu_K(x))  \frac{p_{\lambda}(\overline{x})|x|^{1-n} G_t(\overline{x}, |x|)}   { \depsi_t(\nu_K(x), x\cdot \nu_K(x))}\,dx  \label{new measue-11-27-K},   \\ 
             \int_{\sphere} g(u)\, dC_{\Theta}(K, u)&=\int_{\partial K}g(\nu_K(x))  \frac{ [x\cdot \nu_K(x)]^{2} p_{\lambda}(\overline{x}) |x|^{1-n} G_t(\overline{x}, |x|)}   { \depsi_t(\nu_K(x), [x\cdot \nu_K(x)]^{-1})}\,dx.\label{nn-K}
  \end{align}

 We now prove a variational formula to derive the  Musielak-Orlicz-Gauss image measure. Let $\lambda$ be a Lebesgue measure on $\sphere$. Suppose that $\wp:   \sphere \times (0, \infty)\to \R$ is a function  such that, for
$K\in\cKo$, the function
$\xi\mapsto\wp(\xi, \rho_K(\xi))$ is measurable on $\sphere$ and is integrable with respect to $\lambda$. Define,  $\VGL(K)$,   the general dual volume of $K\in \cKo$ with respect to $\lambda$,  by  \begin{align}\label{vlam}
        \VGL(K)=\int_{\sphere} \wp(\xi, \rho_K(\xi))\, d\lambda(\xi).
        \end{align} In general, one can define $\VGL$ for  all $f\in  C^+(\sphere)$, where $C^+(\Omega)$ for $\Omega\subseteq \sphere$ denotes the set of all positive
continuous functions defined on $\Omega$. That is, if $f\in  C^+(\sphere)$, 
\begin{align}\label{vlam2}
    \VGL(f)=\int_{S^{n-1}}G(\xi, f(\xi))\,d\lambda(\xi).
    \end{align} Clearly, $\widetilde{V}_{G,\lambda}(K)=\widetilde{V}_{G,\lambda}(\rho_K)$ for
$K\in \cKo$. When $\,d\lambda=\,d\xi$, $\VGL$ becomes the general dual volume $\dvetV(K)$ \cite{GHWXY, GHXY} given by  \begin{align*} 
       \dvetV(K)=\int_{\sphere} \wp(\xi, \rho_K(\xi))\, d\xi.
        \end{align*}

Hereafter, $\Omega\subseteq \sphere$ is always assumed to be a
closed set not contained in any closed hemisphere of $\sphere$.  Denote by $C(\Omega)$ the set of all continuous functions defined on
$\Omega\subseteq \sphere$. For each $f\in C^+(\Omega)$, one can define two convex bodies associated to $f$:   the
Wulff shape generated by  $f$  \begin{align}\label{july287}
[f]= \bigcap_{\xi\in\Omega}\big\{x\in\R^n: x \cdot \xi \leq f(\xi)\big\},
\end{align} and the convex hull generated by  $f$:  \begin{align}\label{conv}
\langle f\rangle =\conv\big\{f(\xi)\xi: \xi\in\Omega\big\}. \end{align} Here, $\conv(E)$ denotes the convex hull of set $E\subseteq \Rn$, i.e., $\conv(E)$ is the smallest convex closed set containing $E$. It is easily checked that $[f]\in \cKo$ and $\langle f\rangle \in \cKo$ for $f\in C^+(\Omega)$.
A fundamental relation between the Wulff shape and the convex hull is \begin{align}\label{relation}
[f]^{*}=\langle 1/f \rangle
\end{align}   for $f\in C^+(\Omega)$ (see e.g., \cite[Lemma 2.8]{LYZActa}).  Obviously, for $f\in C^+(\Omega)$,
\begin{align}\label{le1}h_{[
f]}\leq f \ \ \mathrm{and}\ \ \rho_{\langle
f\rangle}\geq f \ \ \mathrm{on} \ \Omega.\end{align} In particular, if $\Omega=\sphere$ and $K\in \cKo$,
\begin{align}\label{hk}[h_K]=K \ \ \mathrm{and}\ \ \langle \rho_{K}\rangle =K .\end{align}

  Let $\depsi\in \cC_I\cup\cC_d$. Recall that 
$\psi_{\xi}=\depsi(\xi, \cdot)$ for each fixed $\xi\in \sphere$.  Then $\psi_{\xi}(\cdot)$ is a strictly monotonic function on $(0,\infty)$, and hence $\psi_{\xi}^{-1}(\cdot)$ exists and is also strictly monotonic.  Let  $f\in C^+(\Omega)$ and $g\in C(\Omega)$. As $\Omega\subseteq \sphere$ is compact and $\depsi_t$ is continuous, there exists a constant $\delta>0$, such that, for all $(\xi,\varepsilon)\in \Omega\times(-\delta,\delta)$,   it is meaningful to define  the {\em Musielak-Orlicz
addition} by
\begin{align}\label{genplus21} f_{\varepsilon}(\xi)= \psi_{\xi}^{-1}\left(\psi_{\xi}(f(\xi))+\varepsilon g(\xi)\right).
\end{align}  
Clearly, $f_{0}(\xi)=f(\xi)$ and $ f_{\varepsilon}(\xi) \in C^+(\Omega)$. Using the chain rule, it is easy
to verify that,   \begin{align}\label{der}
    \frac{\partial f_{\varepsilon}(\xi)
    }{\partial \varepsilon}\Big|_{\varepsilon=\varepsilon_0}
    = \frac{g(\xi)}{\psi'_{\xi} (f_{\varepsilon_0}(\xi))}
\end{align}  for any $(\xi,\varepsilon_0)\in
\Omega\times(-\delta,\delta)$, where $$\psi'_{\xi}(t)=\depsi_t(\xi, t)=\frac{\partial \depsi(\xi, t)}{\partial t}.$$
Moreover, due to the compactness of  $\Omega$, a standard argument shows that $f_{\varepsilon}\to f$  
uniformly on $\Omega$ as $ \varepsilon\rightarrow 0$.

We shall need the following lemma.  
 
\bl \label{variation--11-21}Let $\Omega\subseteq \sphere$ be a closed
set that is not contained in any closed hemisphere of $\sphere$. Let   $f_{\varepsilon}$ be given as in \eqref{genplus21} with  $f\in C^+(\Omega)$, $g\in C(\Omega)$ and
$\depsi \in \cC_I\cup\cC_d$. 

\vskip 2mm \noindent i) 
For $v\in
S^{n-1}\setminus \eta_{\langle f\rangle}$, let $u_0=\alpha_{\langle f\rangle^*}(v)$ and then 
\begin{align}\label{hf}
\lim_{\varepsilon\rightarrow 0}\frac{\log h_{\langle
f_{\varepsilon}\rangle}(v)-\log h_{\langle
f\rangle}(v)}{\varepsilon}=\frac{g(u_0)}{f(u_0)\,
\depsi_t(u_0, f(u_0))}.
\end{align}  ii)  For $\xi\in S^{n-1}\setminus \eta_{\langle 1/f\rangle}$, let $u_1=\alpha_{[f]}(\xi)$ and then
\begin{align}\label{414}
\lim_{\varepsilon\rightarrow 0}\frac{\log
\rho_{[f_{\varepsilon}]}(\xi)-\log
\rho_{[f]}(\xi)}{\varepsilon}=\frac{g(u_1)}{f(u_1)\,
\depsi_t(u_1, f(u_1))}.
\end{align} 
\el 

\begin{proof}  From \eqref{der}, for any $(\xi,\varepsilon_0)\in \Omega\times(-\delta,\delta)
$, we have
\begin{align*}
    \frac{\partial \log f_{\varepsilon}(\xi)
    }{\partial \varepsilon}\Big|_{\varepsilon=\varepsilon_0}
    =
    \frac{g(\xi)}{f_{\varepsilon_0}(\xi)\psi_{\xi}^{\prime}(f_{\varepsilon_0}(\xi))}.
\end{align*}
By the mean value theorem, there exists $\theta(\xi,\varepsilon)\in
(0,1)$ such that
\begin{align}\label{hope1}
\log f_{\varepsilon}(\xi)-\log f(\xi)=\varepsilon\,\frac{g(\xi)}
{f_{\theta(\xi,\varepsilon)\varepsilon}(\xi)\,\psi_{\xi}^{\prime}(f_{\theta(\xi,\varepsilon)\varepsilon}(\xi))}.
\end{align}

Recall that 
$\alpha_{\langle f\rangle^*}(S^{n-1}\setminus \eta_{\langle
f\rangle})\subseteq\Omega$ \cite[(4.24)]{LYZActa}.  Let $v\in\sphere\setminus \eta_{\langle f\rangle}$. Then there exists a vector, say $v_0\in \sphere$, such that, for $\varepsilon\in (-\delta, \delta)$, \begin{align}\label{u0}
h_{\langle
f\rangle}(v)= ( v_{0}\cdot v) f(v_{0}) \quad \textrm{and} \quad  h_{\langle f_{\varepsilon}\rangle}(v)\geq ( v_{0}\cdot v)
f_{\varepsilon}(v_{0}).
\end{align}  From \eqref{star} and  \eqref{u0}, one clearly has $v_0=\alpha^*_{\langle f\rangle}(v)=\alpha_{\langle f\rangle^*}(v)$ since $f(v_0)v_0\in \partial \langle f\rangle$ and $v\in\sphere\setminus
\eta_{\langle f\rangle}$.  As $\alpha_{\langle f\rangle^*}$  is injective on $\sphere\setminus
\eta_{\langle f\rangle}$, one further has $v_0=u_0\in \Omega$.  
 It follows from \eqref{hope1} and \eqref{u0} that
\begin{align}
\log h_{\langle f_{\varepsilon}\rangle}(v)-\log h_{\langle
f\rangle}(v)
\geq \log f_{\varepsilon}(u_{0})-\log f(u_{0}) =\varepsilon\,\frac{g(u_{0})}
{f_{\theta(u_{0},\varepsilon)\varepsilon}(u_{0})\,\psi_{u_{0}}^{\prime}(f_{\theta(u_{0},\varepsilon)\varepsilon}(u_{0}))}.
\label{comparison-11-21-1}
\end{align}  By \eqref{conv},  there exists a $u_{\varepsilon}\in\Omega$ such that
\begin{align}\label{rel}
h_{\langle f_{\varepsilon}\rangle}(v)= (u_{\varepsilon}\cdot v)
f_{\varepsilon}(u_{\varepsilon})\quad \textrm{and} \quad h_{\langle
f\rangle}(v)\geq (u_{\varepsilon}\cdot v) f(u_{\varepsilon}).
\end{align}
Thus, from   \eqref{hope1} and \eqref{rel}, we have
\begin{eqnarray}\label{comparison-11-21}
\log h_{\langle f_{\varepsilon}\rangle}(v)-\log h_{\langle
f\rangle}(v)
\leq\log f_{\varepsilon}(u_{\varepsilon})-\log  f(u_{\varepsilon}) 
=\varepsilon\,\frac{g(u_{\varepsilon})}
{f_{\theta(u_{\varepsilon},\varepsilon)\varepsilon}(u_{\varepsilon})\,\psi_{u_{\varepsilon}}^{\prime}(f_{\theta(u_{\varepsilon},\varepsilon)\varepsilon}(u_{\varepsilon}))}.
\end{eqnarray} Note that $\lim_{\varepsilon\rightarrow0} u_{\varepsilon}=u_{0}$ which is a direct consequence of the fact that $f_{\varepsilon}\rightarrow f$ uniformly on $\Omega$; this can be easily proved following along the same lines as the proof of formula (4.8) in \cite{LYZActa}.  
 By letting $\varepsilon\to 0$, the desired argument  \eqref{hf} follows immediately from \eqref{comparison-11-21-1}, \eqref{comparison-11-21}, 
 the continuity of $g$ and $\depsi_t,$ and $\theta(\cdot, \varepsilon)\in (0, 1)$ for all $\varepsilon\in (-\delta, \delta)$.
 
 Now let us prove \eqref{414}. 
For any $\xi\in\sphere\setminus \eta_{\langle 1/f\rangle}$ and
$\varepsilon\in (-\delta, \delta)$, it follows from \eqref{bi-polar--12} and  \eqref{relation} that
\begin{align}\label{f}
\rho_{[f_{\varepsilon}]}(\xi)=\rho_{\langle
\frac{1}{f_{\varepsilon}}\rangle^*}(\xi)=\frac{1}{h_{\langle
\frac{1}{f_{\varepsilon}}\rangle}(\xi)}.
\end{align}
Recall that $\widetilde{\depsi}(\xi, t)=\depsi(\xi, 1/t)$ and let
$\widetilde{\psi}_{\xi}(t)=\psi_{\xi}(1/t)$ for $t\in (0, \infty)$. Clearly, $\widetilde{\psi}_{\xi}(\cdot)$ is a monotonic function on $(0, \infty)$ if $\depsi\in \cC_I\cap\cC_d$. 
 Hence,   it is meaningful to rewrite 
\eqref{genplus21} as follows: 
\begin{align}\label{til}
\frac{1}{f_{\varepsilon}(\xi)}=\widetilde{\psi}_{\xi}^{-1}\left(\widetilde{\psi}_{\xi}\left(\frac{1}{f(\xi)}\right)+\varepsilon
g(\xi)\right). 
\end{align} It is easily checked that
$\widetilde{\psi}_{\xi}'(t)=-t^{-2}\psi_{\xi}'(1/t)$.  By  \eqref{relation}, \eqref{hf} (in fact, with $u_0$ replaced by $\alpha_{\langle 1/f\rangle^*}(\xi)=\alpha_{[f]}(\xi)=u_1$ due to  \eqref{relation}), \eqref{f} and \eqref{til},  one has  
\begin{align*}
\lim_{\varepsilon\rightarrow 0}\frac{\log
\rho_{[f_{\varepsilon}]}(\xi)-\log \rho_{[f]}(\xi)}{\varepsilon} =
-\lim_{\varepsilon\rightarrow 0}\frac{\log h_{\langle
1/f_{\varepsilon}\rangle}(\xi)-\log h_{\langle
1/f\rangle}(\xi)}{\varepsilon} =\frac{g(u_1)}{f(u_1)\,
\depsi_t(u_1, f(u_1))}.
\end{align*}
This completes the proof of \eqref{414}.  \end{proof}

Note that $[f_i]\rightarrow [f]$ and $\langle f_i\rangle \rightarrow \langle f\rangle$,  if   $f\in C^+(\Omega)$ and $f_i\in C^+(\Omega)$  for all $i\in \N$ such that $f_i\to f$ uniformly on $\Omega$ (see e.g., \cite [p.345]{LYZActa} and  \cite[Lemma 7.5.2]{Sch}). It is also true that $K_i\rightarrow K$ with  $K_i\in \cKo$ for all $i\in \N$ and $K\in \cKo$ is equivalent to $\rho_{K_i}\rightarrow \rho_K$ uniformly on $\sphere$. An application of the dominated convergence theorem yields the continuity of $\VGL(\langle f_{\varepsilon}\rangle^*)$ and $\widetilde{V}_{G,\lambda}([f_{\varepsilon}])$ on $\varepsilon\in (-\delta, \delta)$. The following theorem provides a result on the variational formulas regarding $\VGL(\langle f_{\varepsilon}\rangle^*)$ and $\widetilde{V}_{G,\lambda}([f_{\varepsilon}])$, which  can be used to derive the Musielak-Orlicz-Gauss image measure and its polar.  

  \bt\label{ovev-cor} Let $\Omega\subseteq \sphere$ be a closed
set that is not contained in any closed hemisphere of $\sphere$.  Let $\Theta=(G, \depsi, \lambda)$ be a triple such that $G\in \cC$,  $\depsi \in \cC_I\cup\cC_d$, and $\lambda\in \cM.$ For  $f_{\varepsilon}$ given  by \eqref{genplus21} with  $f\in C^+(\Omega)$ and $g\in C(\Omega)$, one has  
 \begin{align}\label{rr}
    \lim_{\varepsilon\rightarrow 0}\frac{\VGL(\langle f_{\varepsilon}\rangle^*)-\VGL(\langle f\rangle^*)}{\varepsilon} &=
    -\int_{\Omega} g(u)\, dC_{\Theta}(\langle f\rangle^*, u), \\ \label{variation-11-27-12}
    \lim_{\varepsilon\rightarrow 0}\frac{\widetilde{V}_{G,\lambda}([f_{\varepsilon}])-\widetilde{V}_{G,\lambda}([f])}{\varepsilon}&=
    \int_{\Omega} g(u)\, d\deV([f], u).
    \end{align}
    \et
  
  \begin{proof}  It has been shown in \cite[p.17]{GHWXY} that for any $h_0\in C^+(\Omega)$, one has  
$$ h_{[h_0]}(\alpha_{[h_0]}(\xi))=h_0(\alpha_{[h_0]}(\xi)) \quad \mathrm{for}\  \cH^{n-1}-\mathrm{almost\ 
all}\ \xi\in \sphere.
$$ Applying this result to $h_0=1/f$ for $f\in C^+(\Omega)$, one can obtain that, by  
\eqref{bi-polar--12}, \eqref{relation} and  the fact that $\lambda$ is absolutely continuous  with respect to $\,d\xi$,  
 $$ \rho_{\langle
f\rangle}(\alpha_{\langle f\rangle^*}(\xi))=f(\alpha_{\langle
f\rangle^*}(\xi))\quad \mathrm{for}\  \lambda-\mathrm{almost\ 
all}\ \xi\in \sphere.
$$   

For $\varepsilon\in (-\delta, \delta)$, let $K_{\varepsilon}=\langle f_{\varepsilon}\rangle^*=[1/f_{\varepsilon}]\in \cKo$ where $f_{\varepsilon}\in C^+(\Omega)$ is given  by \eqref{genplus21} with  $f\in C^+(\Omega)$ and $g\in C(\Omega)$. When $\varepsilon=0$, $K_{0}=\langle f\rangle^*=[1/f]\in \cKo$. By  
\eqref{bi-polar--12}, one sees that $(K^*)^*=K$ for $K\in \cKo$ and hence  for all $\varepsilon\in (-\delta, \delta)$,  $$\rho_{K_{\varepsilon}}=\frac{1}{h_{K_{\varepsilon}^*}}=\frac{1}{h_{\langle f_{\varepsilon}\rangle}}.$$ This further implies  that  for each $\xi\in \sphere$, 
\begin{align}
\frac{\partial G(\xi, \rho_{K_{\varepsilon}}(\xi))}{\partial  \varepsilon} &= \rho_{K_{\varepsilon}}(\xi)
G_t(\xi, \rho_{K_{\varepsilon}}(\xi)) \frac{\,d \log \rho_{K_{\varepsilon}}(\xi)}{\,d\varepsilon}  \nonumber \\  &=- \left(\frac{1}{h_{\langle f_{\varepsilon}\rangle}(\xi)}  
G_t\Big(\xi, \frac{1}{h_{\langle f_{\varepsilon}\rangle}(\xi)} \Big) \frac{\,d \log h_{\langle f_{\varepsilon}\rangle}(\xi) }{\,d\varepsilon}\right).
\label{chain}
\end{align} Hence, for $\lambda$-almost all $\xi\in \sphere$, by \eqref{hf} and \eqref{chain}, one can get  
\begin{align} \frac{\partial G(\xi, \rho_{K_{\varepsilon}}(\xi))}{\partial \varepsilon} \bigg|_{\varepsilon=0}  &=- \left(\frac{1}{h_{\langle f\rangle}(\xi)}  
G_t\Big(\xi, \frac{1}{h_{\langle f\rangle}(\xi)} \Big) \frac{g(\alpha_{\langle f\rangle^*}(\xi))}{f(\alpha_{\langle f\rangle^*}(\xi))\,
\depsi_t(\alpha_{\langle f\rangle^*}(\xi), f(\alpha_{\langle f\rangle^*}(\xi)))}\right) \nonumber \\ &=-   \left(\frac{\rho_{K_0}(\xi)  
G_t\big(\xi,  \rho_{K_0}(\xi) \big) g(\alpha_{K_0}(\xi))}{f(\alpha_{K_0}(\xi))\,
\depsi_t(\alpha_{K_0}(\xi), f(\alpha_{K_0}(\xi)))}\right). \label{chain-7-31-1}
\end{align} Moreover, there are $\delta_0\in (0,\delta)$ and a constant $M>0$,
such that for $\varepsilon\in (-\delta_0,\delta_0)$ and all
$\xi\in\sphere$, 
\begin{align*}
\left|\frac{G(\xi, \rho_{K_{\varepsilon}}(\xi))-G(\xi, \rho_{K_{0}}(\xi))}{\varepsilon}\right|\leq M.
\end{align*}   This can be done because, involving  in \eqref{comparison-11-21-1}, \eqref{comparison-11-21} and \eqref{chain}, the sets (i.e., $\Omega$ and $\sphere$) are all compact, the functions $G\in \cC$ and  $\depsi\in \cC_I\cup \cC_d$ (hence  $G_t$ and $\depsi_t$ are all continuous),  
and the family of functions $f_{\varepsilon}$ is uniformly bounded on $\Omega$ from below by a positive number and from above by a finite number. For more details on how to find $M$, please refer to, e.g.,  the proofs of  Lemma 4.2 in \cite{LYZActa},  Theorem 4.1 \cite{XY2017-1} and  Lemma 4.2 in \cite{ZSY2017}. It follows from \eqref{chain-7-31-1} and the  dominated convergence theorem that, for $K_{\varepsilon}=\langle f_{\varepsilon}\rangle^*$ and $K_{0}=\langle f\rangle^*$,   
\begin{align}
\lim_{\varepsilon\rightarrow
0}\frac{\widetilde{V}_{G,\lambda}(\langle
f_{\varepsilon}\rangle^*)-\widetilde{V}_{G,\lambda}(\langle
f\rangle^*)}{\varepsilon}  &=\lim_{\varepsilon\rightarrow
0}\frac{\widetilde{V}_{G,\lambda}(K_{\varepsilon})-\widetilde{V}_{G,\lambda}(K_{0})}{\varepsilon}  \nonumber \\ & =
 \lim_{\varepsilon\rightarrow 0}  \int_{\sphere}\frac{G(\xi, \rho_{K_{\varepsilon}}(\xi))-G(\xi, \rho_{K_{0}}(\xi))}{\varepsilon}
 \,d\lambda(\xi)
 \nonumber\\&=  \int_{\sphere} \lim_{\varepsilon\rightarrow 0}  \frac{G(\xi, \rho_{K_{\varepsilon}}(\xi))-G(\xi, \rho_{K_{0}}(\xi))}{\varepsilon}
 \,d\lambda(\xi)
 \nonumber\\
&=-\int_{\sphere\setminus \eta_{\langle f\rangle }}\frac{\rho_{K_0}(\xi)  
G_t\big(\xi,  \rho_{K_0}(\xi) \big) g(\alpha_{K_0}(\xi))}{f(\alpha_{K_0}(\xi))\,
\depsi_t(\alpha_{K_0}(\xi), f(\alpha_{K_0}(\xi)))}\,d\lambda(\xi). \label{var-7-31-2-1}
\end{align}

 Recall that 
$\alpha_{\langle f\rangle^*}(S^{n-1}\setminus \eta_{\langle
f\rangle})\subseteq\Omega$ \cite[(4.24)]{LYZActa}.  
The compactness of $\Omega$, together with the Tietze extension theorem, yields the existence of  $\overline{g}: S^{n-1}\rightarrow \R$ such
that $\overline{g}$ is continuous on $\sphere$ and   for $\xi\in S^{n-1}\setminus \eta_{\langle f\rangle}$,
$$g(\alpha_{\langle f\rangle^*}(\xi))= g(\alpha_{K_0}(\xi))= (\overline{g}1_{\Omega})(\alpha_{K_0}(\xi))=(\overline{g}1_{\Omega})(\alpha_{\langle f\rangle^*}(\xi)),
$$ where $1_{E}$ is the indicator function of $E$, i.e., $1_{E}(x)=1$ if $x\in E$ and $1_{E}(x)=0$ if $x\notin E$.  Applying this to \eqref{var-7-31-2-1}, one further gets  
\begin{align}
\lim_{\varepsilon\rightarrow
0}\frac{\widetilde{V}_{G,\lambda}(\langle
f_{\varepsilon}\rangle^*)-\widetilde{V}_{G,\lambda}(\langle
f\rangle^*)}{\varepsilon}    
&=-\int_{\sphere} \frac{(\overline{g}1_{\Omega})(\alpha_{K_0}(\xi)) \rho_{K_0}(\xi)  
G_t\big(\xi,  \rho_{K_0}(\xi) \big) }{f(\alpha_{K_0}(\xi))\,
\depsi_t(\alpha_{K_0}(\xi), f(\alpha_{K_0}(\xi)))} \,d\lambda(\xi)\nonumber\\
&=-\int_{\sphere}(\overline{g}1_{\Omega})(u) \,
dC_{\Theta}(K_0,
u)\nonumber\nonumber\\
&=-\int_{\Omega}g(u) \, dC_{\Theta}(\langle
f\rangle^*, u).\nonumber
\end{align} where we have used \eqref{nn} in the second equality and   $$f(\alpha_{K_0}(\xi))=f(\alpha_{\langle f\rangle^*}(\xi))=\rho_{\langle f\rangle}(\alpha_{\langle f\rangle^*}(\xi))=\rho_{K_0^*}(\alpha_{K_0}(\xi)).$$ 
This concludes the proof of \eqref{rr}.

The variational formula \eqref{variation-11-27-12} follows along the same lines as the proof for \eqref{rr}, based on  \eqref{414}. A more
direct proof for \eqref{variation-11-27-12} can  be given by
the combination of \eqref{wc} and \eqref{rr}. Indeed, let $\Theta=(G, \depsi, \lambda)$ be a given triple and  $\widetilde{\Theta}=(G, \widetilde{\depsi}, \lambda)$. It follows from \eqref{wc}, 
\eqref{relation}  and \eqref{rr}  (applied to  $\widetilde{\Theta}$ instead of $\Theta$ due to \eqref{til})  that
\begin{align}
 \lim_{\varepsilon\rightarrow
0}\frac{\widetilde{V}_{G,\lambda}([f_{\varepsilon}])-\widetilde{V}_{G,\lambda}([f])}{\varepsilon}
 &= \lim_{\varepsilon\rightarrow
0}\frac{\widetilde{V}_{G,\lambda}(\langle
1/f_{\varepsilon}\rangle^*)-\widetilde{V}_{G,\lambda}(\langle
1/f\rangle^*)}{\varepsilon}
 \nonumber\\&=  -\int_{\Omega} g(u)\, dC_{\widetilde{\Theta}}(\langle 1/f\rangle^*, u)
   \nonumber\\&=  \int_{\Omega} g(u)\, d\widetilde{C}_{\Theta}(\langle 1/f\rangle^*, u)
    \nonumber\\&=  \int_{\Omega} g(u)\, d\widetilde{C}_{\Theta}([f], u).\nonumber
\end{align} This completes the proof of \eqref{variation-11-27-12}.\end{proof}

     \section{The Musielak-Orlicz-Gauss image problem}\label{M:4} \setcounter{equation}{0}
   This section is dedicated to introduce the  Musielak-Orlicz-Gauss image problem and some of its special cases. Indeed, it has been explained in Section \ref{Section-var} that both $\deV(K,
\cdot)$ and $C_{\Theta}(K,\cdot)$ are signed Borel measures on $\sphere$.  We now prove some basic properties for  $C_{\Theta}(K,\cdot)$ and $\deV(K,\cdot)$ for $K\in \cKo$. 
  
   \bp\label{prop 11-28--1}  Let $\Theta=(G, \depsi, \lambda)$ be a triple such that $G\in \cC$,  $\depsi \in \cC_I\cup\cC_d$ and $\lambda\in \cM.$   Let $K\in \cKo$.  Then the following statements hold.               \vskip 1mm \noindent  i)  Both  $C_{\Theta}(K_i,\cdot)\rightarrow C_{\Theta}(K,\cdot)$ and $\deV(K_i,\cdot)\rightarrow \deV(K,\cdot)$ weakly for any sequence of $\{K_i\}_{i\in \N}$ such that $K_i\in \cKo$ for any $i\in \N$ and $K_i\rightarrow K\in \cKo$.   
             \vskip 1mm \noindent  ii) The signed measures $C_{\Theta}(K,\cdot)$ and $\deV(K,\cdot)$ are absolutely continuous with respect to $S(K,\cdot)$.
 \vskip 1mm \noindent  iii)  If  $G$ and $\depsi$ are either both in $\cC_I$ or both in $\cC_d$, then $\deV(K,\cdot)$
 and $C_{\Theta}(K,\cdot)$  are nonzero finite Borel
measures.  If, in addition, $\lambda$ is strictly positive on nonempty open subsets of $\sphere$, then $\deV(K,\cdot)$
 and $C_{\Theta}(K,\cdot)$  are not concentrated on any closed hemisphere of
$\sphere$.  The same arguments also hold for $-\deV(K,\cdot)$ and  $-C_{\Theta}(K,\cdot)$, if one of $G$ and $\depsi$ is in $\cC_I$ and the other one is in $\cC_d$.

  \ep \begin{proof} Due to \eqref{wc},   only $\deV(K,
\cdot)$ will be discussed.   Part i) follows easily from a standard argument of the dominated convergence theorem, by \eqref{relation-G} and the facts that $\rho_{K_i}\to \rho_K$ and $h_{K_i}\to h_K$ uniformly on $\sphere$,  $\alpha_{K_i}\rightarrow \alpha_K$  (holding except a subset of $\sphere$ whose  $\cH^{n-1}$-measure is zero \cite[Lemma 2.2]{LYZActa}),  and  $\lambda(K, \cdot)$ (and hence $\lambda^*(K, \cdot))$ is weakly convergent on $\cKo$ (\cite[Lemma 3.4]{BLYZZ2020}).   

\vskip 2mm \noindent ii) Let $K\in \cKo$ and $\omega\in \cB$ such that $S(K, \omega)=0$. It has been proved in \cite[Lemma 3.5]{BLYZZ2020} that $\lambda(K, \cdot)$ for $K\in \cKo$ is absolutely continuous with respect to the surface area measure $S(K^*, \cdot)$.  Applying this to $\polar$, one gets  $\lambda(\polar, \omega)=\lambda^*(K, \omega)=0$.  As $K\in \cKo$,  there exist two constants $0<r_0<R_0<\infty$ such that  both $h_{K}$ and $\rho_{K}$ are in $(r_0, R_0)$ on $\sphere$. The continuity of
$G_{t}$ and $\depsi_{t}$ yield that
\begin{align}c_1:=\sup_{\xi\in \sphere} \left|  \frac{\rho_K(\xi) G_t(\xi, \rho_K(\xi))}{h_{K}(\alpha_K(\xi))\depsi_t(\alpha_K(\xi), h_{K}(\alpha_K(\xi)))}\right|<\infty. \label{upp-7-31-1}\end{align} Together with  \eqref{gencdef-7-1-115} (or \eqref{new measue-11-27}) and \eqref{gencdef-7-1-115-uu}, one has 
         \begin{align} \big| \deV(K, \omega) \big| &= \bigg| \int_{\pmb{\alpha}^{*}_K(\omega)} \frac{\rho_K(\xi) G_t(\xi, \rho_K(\xi))}{h_{K}(\alpha_K(\xi))\depsi_t(\alpha_K(\xi), h_{K}(\alpha_K(\xi)))}\,d\lambda(\xi)\bigg| \nonumber \\
                    &\leq c_1\int_{{\pmb{\alpha}^{*}_K(\omega)}}\,d\lambda(\xi)= c_1\lambda\left(\boldsymbol{\alpha}^{*}_{K}(\omega)\right)=c_1\lambda^*(K, \omega)=0. \label{non-conc-2}
                \end{align} This concludes that  $\deV(K,\cdot)$ is absolutely continuous with respect to $S(K,\cdot)$.

\vskip 2mm \noindent iii) Only the proof for the case when $G\in \cC_I$ and $\depsi\in \cC_I$  will be given, and the other cases follow along the same lines.  A calculation similar to \eqref{upp-7-31-1} yields that 
\begin{align}\label{small}
c_2:=\inf_{\xi\in\sphere} \frac{\rho_K(\xi) G_t(\xi, \rho_K(\xi))}{h_{K}(\alpha_K(\xi))\depsi_t(\alpha_K(\xi), h_{K}(\alpha_K(\xi)))}>0.
\end{align} This implies that $\deV(K, \cdot)$ is a nonzero measure. 
Following the proof of \eqref{non-conc-2}, one can also prove  \begin{align*} \deV(K, \sphere)  \leq c_1\lambda^*(K, \sphere)<\infty.\end{align*} Hence $\deV(K, \cdot)$ is  finite.

Assume that, in addition, $\lambda$ is strictly positive on nonempty open subsets of $\sphere$.  We now claim that $\deV(K,\cdot)$ satisfies  \eqref{not-concentration-1}. This is an easy consequence of \eqref{new measue-11-27}, \eqref{small} and Lemma
\ref{conc}:   for any $u\in\sphere$,
                \begin{align*}
                    \int_{\sphere}( u\cdot v)_+\, d\deV(K, v)&=
                    \int_{\sphere}( u \cdot\alpha_K(\xi))_+ \ \frac{\rho_K(\xi) G_t(\xi, \rho_K(\xi))}{h_{K}(\alpha_K(\xi))\depsi_t(\alpha_K(\xi), h_{K}(\alpha_K(\xi)))}\,d\lambda(\xi)\nonumber\\
&\geq c_2 \int_{\sphere}(u \cdot\alpha_K(\xi))_+ \,d\lambda(\xi)\nonumber>0, 
                \end{align*}
where the second inequality follows from \eqref{star} and \eqref{not-concentration-1-22}  (applying to $\polar$).    \end{proof}

  Proposition \ref{prop 11-28--1} suggests the following Musielak-Orlicz-Gauss image problem.  
      
\begin{problem}[The Musielak-Orlicz-Gauss image problem] \label{MOGIP} Let  $G\in \cC$, $\depsi\in \cC$, and  
$\lambda$ be a nonzero finite Lebesgue measure on $\sphere$.  Under what conditions on $\Theta=(G, \Psi, \lambda)$ and a nonzero finite Borel measure $\mu$ on $\sphere$ do there exist a $K\in \cKo$ and a constant $\tau\in \R$ such that  $\mu =\tau \deV(K,\cdot)?$
\end{problem} Let $|\mu|=\int_{\sphere}\,d\mu(u)$. 
   Clearly, if Problem \ref{MOGIP} has $K\in \cKo$ as its solution, then $$\tau=\frac{|\mu|}{\deV(K,\sphere)}.$$ 
Problem \ref{MOGIP} is for the characterization of the  Musielak-Orlicz-Gauss image measure. Similar problem can be posed for the  polar Musielak-Orlicz-Gauss image measure. 
  \begin{problem}[The  polar Musielak-Orlicz-Gauss image problem] \label{MOGIP-p} Let  $G\in \cC$, $\depsi\in \cC$, and  
$\lambda$ be a nonzero finite Lebesgue measure on $\sphere$.  Under what conditions on $\Theta=(G, \Psi, \lambda)$  and a nonzero finite Borel measure $\mu$ on $\sphere$ do there exist a $K\in \cKo$ and a constant $\kappa\in \R$ such that $\mu =\kappa C_{\Theta} (K,\cdot)?$
\end{problem}  Again if Problem \ref{MOGIP-p} has $K\in \cKo$ as its solution, then $$\kappa=\frac{|\mu|}{C_{\Theta}(K,\sphere)}.$$
 When $G=\depsi=\log t$,  it can be seen from \eqref{relation-G} that  Problems   \ref{MOGIP} and \ref{MOGIP-p} reduce to the Gauss image problem (i.e., Problem \ref{Gauss-I-p}) introduced in \cite{BLYZZ2020}.  From the discussion in Section \ref{Section-var}, one clearly sees that 
 Problem \ref{MOGIP} also generalizes the Minkowski problem \cite{min1897, min1903},   the $L_p$ Minkowski problem \cite{Lu93}, the Orlicz-Minkowski problem  \cite{HLYZ2010}, the ($L_{p}$)  dual Minkowski  \cite{LYZActa, LYZ-Lp}, the  dual Orlicz-Minkowski problems \cite{GHWXY, GHXY, XY2017-1, ZSY2017}, and  the ($L_p$ and Orlicz) Aleksandrov problem \cite{Alexs1942, FH, HLYZ}.

 Problems \ref{MOGIP} and  \ref{MOGIP-p} have close connections with the Monge-Amp\`{e}re  type equations. To see this, let $K\in \cKo$ be smooth enough, in particular, satisfying that $h_K$ is differentiable at each point on $\sphere$ and $\partial K$ has positive Gauss curvature at each point. Denote by  $\nabla h_K(u)$  the gradient of $h_K$ at $u\in \sphere$ and by $\bar{\nabla}h$  the gradient of $h$ with respect to an orthonormal frame on $\sphere$. Then $\nabla h_K=\bar{\nabla}h_K+h_K\iota,$ where  $\iota$ denotes the identity map on $\sphere$ (see, e.g., \cite[(2.2)]{LYZ-Lp}).   Let $\bar{\nabla}^2h$  be the Hessian matrix of $h$ with respect to an orthonormal frame on $\sphere$. Then,  see e.g., \cite[(3.28)]{LYZ-Lp}, for all $u\in \sphere$, $$\frac{\,dS(K, u)}{\,du}=\det (\bar{\nabla}^2h_K(u)+h_K(u)I), $$  where $I$ denotes  the identity matrix. 
 Recall that $\nabla h_K(u)=\nu^{-1}_K(u)$ and $\nabla h_K(\nu_K(x))=x$ hold for all $u\in \sphere$ and $x\in \partial K$. Consequently, by \eqref{surface-8-7}, \eqref{new measue-11-27-K} and \eqref{nn-K}, one gets  \begin{eqnarray}  \, d\deV(K, u) \!\! &=&\!\!   \frac{p_{\lambda}\Big(\frac{\nabla h_K(u)}{|\nabla h_K(u)|}\Big)|\nabla h_K(u)|^{1-n} G_t\Big(\frac{\nabla h_K(u)}{|\nabla h_K(u)|}, |\nabla h_K(u)|\Big)}   { \depsi_t(u,  h_K(u))}\,dS(K, u),   \label{change-v-8-7-1} \\  dC_{\Theta}(K, u)\!\!&=&\!\!    \frac{p_{\lambda}\Big(\frac{\nabla h_K(u)}{|\nabla h_K(u)|}\Big) \big(h_K(u)\big)^2 |\nabla h_K(u)|^{1-n} G_t\Big(\frac{\nabla h_K(u)}{|\nabla h_K(u)|}, |\nabla h_K(u)|\Big)}   { \depsi_t(u,  (h_K(u))^{-1})}\,dS(K, u),
\nonumber \end{eqnarray}  where $\,d\lambda(\xi)=p_{\lambda}(\xi)\,d\xi$ with $p_{\lambda}: \sphere\rightarrow [0, \infty)$. Subsequently, if $\mu$ has its density function to be $p_{\mu}$ with respect to $\,d\xi$,  then \eqref{change-v-8-7-1} yields the following rephrase of Problem \ref{MOGIP} as an Monge-Amp\`{e}re type equation:  \begin{align}\label{new2}
           p_{\mu}= \tau  \frac{P(\bar{\nabla}h+h\iota) \,\det(\bar{\nabla}^2h+hI)}{\depsi_{t}(\cdot, h)}  p_{\lambda}\! \left(\!\frac{\bar{\nabla} h+h \iota}{|\bar{\nabla} h+h \iota|}\!\right), \end{align}
             where $P(y)=|y|^{1-n}G_{t}(\bar{y}, |y|)$ for $y\in \R^n$. Thus,  finding a solution to Problem \ref{MOGIP} requires to find a 
  $\tau\in \R$ and $h:\sphere\to (0,\infty)$ satisfying \eqref{new2}.  Similarly, Problem \ref{MOGIP-p} can be rephrased as follows: 
\begin{align*}
           p_{\mu}= \kappa \frac{h^2 P(\bar{\nabla}h+h\iota) \,\det(\bar{\nabla}^2h+hI)}{\depsi_{t}(\cdot, 1/h)}  p_{\lambda}\! \left(\!\frac{\bar{\nabla} h+h \iota}{|\bar{\nabla} h+h \iota|}\!\right).  \end{align*}

  Some special cases of Problems   \ref{MOGIP} and \ref{MOGIP-p} are of particular interest. 
  
  \vskip 2mm \noindent {\bf Case 1:} $G=t^n/n$, $\Psi\in \cC$, and $\,d\lambda(\xi)=\,d\xi$. In this case, the measure $\deV(K, \cdot)$ will be denoted by $S_{\depsi}(K, \cdot)$ and called the Musielak-Orlicz surface area measure of $K\in \cKo$. Indeed, $S_{\depsi}(K, \cdot)$ has the ($L_p$ and  Orlicz) surface area measures  \eqref{surface-8-7} as its special cases, and satisfies the following formula due to \eqref{gencdef-7-1-115}:
    $$
            S_{\depsi}(K, \omega)=\int_{\pmb{\alpha}^{*}_K(\omega)} \frac{\rho_K^n(\xi)}{h_{K}(\alpha_K(\xi))\depsi_t(\alpha_K(\xi), h_{K}(\alpha_K(\xi)))}\,d\xi. 
         $$ Moreover, letting $u=\alpha_K(\xi)$, it can be verified by \eqref{surface-8-7} and \eqref{vol-change-8-8} that $$\frac{\,dS_{\depsi}(K, u)}{\,dS(K, u)}=\frac{1}{\depsi_t(u, h_K(u))}.$$   A direct consequence of Theorem \ref{ovev-cor} is the following result, which provides a variational formula to derive $S_{\depsi}(K, \cdot)$. 
   \bt\label{ovev-cor-vol} Let   $\depsi \in \cC_I\cup\cC_d$ and $\Omega\subseteq \sphere$ be a closed
set that is not contained in any closed hemisphere of $\sphere$.  For  $f_{\varepsilon}$ given  by \eqref{genplus21} with  $f\in C^+(\Omega)$ and $g\in C(\Omega)$, one has  
$$
    \lim_{\varepsilon\rightarrow 0}\frac{V([f_{\varepsilon}])-V([f])}{\varepsilon}=
    \int_{\Omega} g(u)\, dS_{\depsi}([f], u).
$$
    \et 
  Thus, the following Musielak-Orlicz-Minkowski problem can be proposed.  \begin{problem}[The Musielak-Orlicz-Minkowski problem] \label{MOMP-vol} Under what conditions on    $\depsi\in \cC$ and a nonzero finite Borel measure $\mu$ on $\sphere$ do there exist a $K\in \cKo$ and a constant $\tau\in \R$ such that   $$\mu =\tau S_{\depsi} (K,\cdot)?$$  
\end{problem}   Problem \ref{MOMP-vol}  is related to ``an increasing function" $G\in \cC_I$, and will be studied in our future work \cite{HXYZ-2}.  The Musielak-Orlicz-Minkowski problem deserves its own special attention as it is the direct extension of the $L_p$ and Orlicz Minkowski  problems and lies in the framework of (the extension of) the Brunn-Minkowski theory of convex bodies. By \eqref{new2}, the Monge-Amp\`{e}re type equation related to Problem \ref{MOMP-vol}  is  \begin{align*}  p_{\mu}= \tau   \frac{\det(\bar{\nabla}^2h+hI)}{\depsi_{t}(\cdot, h)}. \end{align*}              
             
   \noindent {\bf Case 2:}  $G\in \cC$, $\Psi\in \cC$, and $\,d\lambda(\xi)=\,d\xi$.  In this case, $\widetilde{V}_{G,\lambda}(K)= \dvetV(K)$ becomes the general dual volume of $K$,  $\widetilde{C}_{G, \depsi}(K, \cdot)$ for $K\in \cKo$ given by \eqref{gencdef-7-1-115--1} defines an Musielak-Orlicz extension of the dual curvature measures \cite{GHWXY, LYZActa, LYZ-Lp, XY2017-1, ZSY2017}, and hence the following dual Musielak-Orlicz-Minkowski problem can be posed. 
  \begin{problem}[The dual Musielak-Orlicz-Minkowski problem] \label{MOMP-dual} Under what conditions on  $G\in \cC$,  $\depsi\in \cC$ and a nonzero finite Borel measure $\mu$ on $\sphere$ do there exist a $K\in \cKo$ and a constant $\tau\in \R$ such that  $\mu =\tau  \widetilde{C}_{G, \depsi}(K,\cdot)?$  
\end{problem} By \eqref{new2}, the corresponding Monge-Amp\`{e}re type equation related to  the dual Musielak-Orlicz-Minkowski problem is   \begin{align*}      p_{\mu}= \tau  \frac{P(\bar{\nabla}h+h\iota) \,\det(\bar{\nabla}^2h+hI)}{\depsi_{t}(\cdot, h)}, \end{align*}
             where $P(y)=|y|^{1-n}G_{t}(\bar{y}, |y|)$ for $y\in \R^n$.

    \vskip 2mm \noindent {\bf  Case 3:} $G=\log t$, $\Psi\in \cC$ and $\lambda$ a nonzero finite Lebesgue measure on $\sphere$. In this case, we shall give the following  Musielak-Orlicz extension of $\lambda^*(K,\cdot)$.
   
     \bd Let $\Theta=(G, \Psi, \lambda)$ be such that   $\lambda$ is a nonzero finite Lebesgue measure on $\sphere$, and  $\depsi\in \cC_I\cup\cC_d$.   For $K\in \cKo$, define $  \widetilde{J}_{\depsi, \lambda}(K, \cdot)=\deV(K, \cdot)$ with $G=\log t$, namely,  for each Borel set $\omega\in \cB$,                $$
           \widetilde{J}_{\depsi, \lambda}(K, \cdot)=\int_{\pmb{\alpha}^{*}_K(\omega)} \frac{1}{h_{K}(\alpha_K(\xi))\depsi_t(\alpha_K(\xi), h_{K}(\alpha_K(\xi)))}\,d\lambda(\xi). 
           $$
        \ed
 Clearly, one can also have,  for all $K\in \cKo$,  \begin{align}\label{equiv-8-25} \widetilde{C}_{(-\log t, \depsi, \lambda)}(K, \cdot)=- \widetilde{J}_{\depsi, \lambda}(K, \cdot).\end{align}  This formula is convenient in later context when finding solutions to Problem \ref{MOMP-log-8-3}.  According to \eqref{relation-G}, it follows that  \begin{align} \frac{\,d \widetilde{J}_{\depsi, \lambda}(K, u)} {\,d\lambda^*(K, u)} = \frac{1}{h_{K}(u)\depsi_t(u, h_{K}(u))}\ \ \ \mathrm{for}\ \ u\in \sphere. \label{relation-int-8-3} \end{align} 
 Clearly, $\widetilde{J}_{\depsi, \lambda}(K, \cdot)$ is a finite signed Borel measure on $\sphere$.  Moreover, it follows from \eqref{new measue-11-27} that,  for any bounded Borel function $g:\sphere\to \R$,    
\begin{align} \int_{\sphere} g(u)\, d\widetilde{J}_{\depsi, \lambda}(K, u) = \int_{\sphere}\frac{g(\alpha_K(\xi))}
            {h_{K}(\alpha_K(\xi))\depsi_t(\alpha_K(\xi), h_{K}(\alpha_K(\xi)))}\,d\lambda(\xi).\label{new measue-11-27-int}
            \end{align} Regarding this measure, one can pose the following problem.  
            
             \begin{problem}\label{MOMP-log-8-3} Let $\depsi\in \cC$ and  $\lambda$ be a nonzero finite Lebesgue measure on $\sphere$.     Under what conditions on  $\Psi,  \lambda$ and a nonzero finite Borel measure $\mu$ on $\sphere$ do there exist a $K\in \cKo$ and a constant $\kappa\in \R$ such that  $\mu =\kappa  \widetilde{J}_{\depsi, \lambda}(K, \cdot)?$
\end{problem} 

  Again, by \eqref{new2}, the corresponding Monge-Amp\`{e}re type equation related to  Problem \ref{MOMP-log-8-3} is   \begin{align*}
           p_{\mu}= \kappa    \frac{ \det(\bar{\nabla}^2h+hI)}{\depsi_{t}(\cdot, h)|\bar{\nabla}h+h\iota|^n }  p_{\lambda}\! \left(\!\frac{\bar{\nabla} h+h \iota}{|\bar{\nabla} h+h \iota|}\!\right), \end{align*}  where $\,d\lambda(\xi)=p_{\lambda}(\xi)\,d\xi$ with $p_{\lambda}: \sphere\rightarrow [0, \infty)$.  Note that $ \widetilde{J}_{\log, \lambda}(K, \cdot)=\lambda^*(K, \cdot)$. Consequently, Problem \ref{MOMP-log-8-3} becomes the Gauss image problem \cite{BLYZZ2020} (up to a difference of polarity of convex bodies). 

A crucial geometric invariant related to $\widetilde{J}_{\depsi, \lambda}(K, \cdot)$ is $\mathcal{E}_{\lambda}(K)$, the entropy of $K\in \cKo$ with respect to the measure $\lambda$. For $K\in \cKo$, let   \begin{align}\label{rela} \mathcal{E}_{\lambda}(K)=
   \widetilde{V}_{\log,\lambda}(K^*)=\int_{S^{n-1}} \log \rho_{K^*}(\xi) d \lambda(\xi).
    \end{align}  Clearly $\mathcal{E}_{\lambda}(\ball)=0$.  When $\,d\lambda(\xi)=\,d\xi$, it reduces to the entropy of $K\in \cKo$, which plays essential roles in solving the Aleksandrov type problems. Letting $G=\log t$ in  \eqref{rr} and \eqref{variation-11-27-12}, by \eqref{new measue-11-27},  \eqref{nn} and \eqref{new measue-11-27-int}, one can easily get the following variational formula.    
    
   \bt\label{ovev-cor-ent} Let $\Omega\subseteq \sphere$ be a closed
set that is not contained in any closed hemisphere of $\sphere$.  Let  $\depsi \in \cC_I\cup\cC_d$ and $\lambda\in \cM$. For  $f_{\varepsilon}$ given  by \eqref{genplus21} with  $f\in C^+(\Omega)$ and $g\in C(\Omega)$, one has  
 \begin{align*}
   \lim _{\varepsilon \rightarrow 0}
\frac{\mathcal{E}_{\lambda}\left(\left\langle
f_{\varepsilon}\right\rangle\right)-\mathcal{E}_{\lambda}\left(\left\langle
f\right\rangle\right)}{\varepsilon}&=
    -\int_{\Omega} g(u)\, d  J_{\depsi, \lambda}  (\langle f\rangle^*, u), \\  
    \lim_{\varepsilon\rightarrow 0}\frac{\mathcal{E}_{\lambda}([f_{\varepsilon}]^*) -\mathcal{E}_{\lambda}([f]^*)}{\varepsilon}&=
    \int_{\Omega} g(u)\, d \widetilde{J}_{\depsi, \lambda}([f], u),
    \end{align*} where $J_{\depsi, \lambda}(K, \cdot)$ is the measure, such that, for all $\omega\in \cB$,   \begin{align}\label{pd-int}
            J_{\depsi, \lambda}(K, \omega)=\int_{\pmb{\alpha}^{*}_K(\omega)} \frac{1}{\rho_{K^*}(\alpha_K(\xi))\depsi_t
            (\alpha_K(\xi), \rho_{K^*}(\alpha_K(\xi)))}\,d\lambda(\xi). 
            \end{align}    \et    Clearly,     $J_{\depsi, \lambda}(K, \cdot)=C_{\Theta}(K, \cdot)$ for $K\in \cKo$ and $\Theta=(\log t, \depsi, \lambda)$.  Formula \eqref{wc} yields that, if $\widetilde{\depsi}(\xi, t)=\depsi(\xi, \frac{1}{t})$ for $(\xi, t)\in \sphere\times (0, \infty)$,   then \begin{align}\label{J-relation}
            J_{\widetilde{\depsi}, \lambda}(K, \cdot )=- \widetilde{J}_{\depsi, \lambda} (K, \cdot).\end{align}  
      
We now state some basic properties for  $ \widetilde{J}_{\depsi, \lambda}(K, \cdot)$ and  $J_{\depsi, \lambda}(K,\cdot)$, which follow from Proposition \ref{prop 11-28--1} by letting $G=\log t$. 
  
   \bp\label{prop 11-28--1-log}  Let $K\in \cKo$, $\depsi \in \cC_I\cup\cC_d$ and $\lambda\in \cM$.  Then the following statements hold.   \vskip 1mm \noindent  i)  Both  $J_{\widetilde{\depsi}, \lambda}(K_i,\cdot)\rightarrow J_{\widetilde{\depsi}, \lambda}(K,\cdot)$ and $\widetilde{J}_{\depsi, \lambda}(K_i,\cdot)\rightarrow \widetilde{J}_{\depsi, \lambda}(K,\cdot)$ weakly for any sequence of $\{K_i\}_{i\in \N}$ such that $K_i\in \cKo$ for any $i\in \N$ and $K_i\rightarrow K\in \cKo$.   
             \vskip 1mm \noindent  ii) Both $J_{\widetilde{\depsi}, \lambda}(K,\cdot)$ and $\widetilde{J}_{\depsi, \lambda}(K,\cdot)$ are absolutely continuous with respect to $S(K,\cdot)$.
 \vskip 1mm \noindent  iii)   If $\depsi\in \cC_I$, then $\widetilde{J}_{\depsi, \lambda}(K, \cdot)$
 and $J_{\widetilde{\depsi}, \lambda}(K,\cdot)$  are nonzero finite Borel
measures.   If, in addition, $\lambda$ is strictly positive on nonempty open subsets of $\sphere$, then $\widetilde{J}_{\depsi, \lambda}(K, \cdot)$
 and $J_{\widetilde{\depsi}, \lambda}(K,\cdot)$ are not concentrated on any closed hemisphere of
$\sphere$.  The same arguments also hold for $-\widetilde{J}_{\depsi, \lambda}(K, \cdot)$ and  $-J_{\widetilde{\depsi}, \lambda}(K,\cdot)$, if  $\depsi\in \cC_d$.   \ep

 \vskip 2mm \noindent {\bf  Case 4:} $G=\log t$,  $\Psi\in \cC$, and $\,d \lambda(\xi)=\,d\xi$.  In this case, Problem \ref{MOGIP}  becomes the Musielak-Orlicz extension of the Aleksandrov problem.  Recall that the Aleksandrov's integral curvature $J(K, \cdot)$ for $K\in \cKo$ is $\lambda(K,\cdot)$ with  $\,d\lambda(\xi)=\,d\xi$. Moreover,  $J^*(K, \cdot)=J(\polar, \cdot)$ for $K\in \cKo$.  Comparing $J(K, \cdot)$ and \eqref{gencdef-7-1-115-uu},  one sees $J^*(K, \cdot)=\widetilde{C}_{\Theta_1}(K, \cdot)$ with $\Theta_1=(\log t, \log t, \,d\xi)$. So $\widetilde{J}_{\depsi}(K, \cdot)=\widetilde{J}_{\depsi, \,d\xi}(K, \cdot)$ defines a Musielak-Orlicz extension of $J^*(K, \cdot)$ and by \eqref{relation-int-8-3},   \begin{align*} \frac{\,d \widetilde{J}_{\depsi}(K, u)} {\,dJ^*(K, u)} = \frac{1}{h_{K}(u)\depsi_t(u, h_{K}(u))}\ \ \ \mathrm{for}\ \ u\in \sphere.\end{align*} Thus, the following Musielak-Orlicz-Aleksandrov problem can be posed; this provides an extension of the Aleksandrov problems \cite{Alexs1942, FH, HLYZ} (again, up to a difference of polarity of convex bodies).    \begin{problem}[The Musielak-Orlicz-Aleksandrov problem] \label{MOMP-Alex}     Under what conditions on  $\Psi$ and a nonzero finite Borel measure $\mu$ on $\sphere$ do there exist a $K\in \cKo$ and a constant $\kappa\in \R$ such that  $\mu =\kappa  \widetilde{J}_{\depsi}(K, \cdot)?$
\end{problem}   Again, by \eqref{new2}, the corresponding Monge-Amp\`{e}re type equation related to the Musielak-Orlicz-Aleksandrov problem is   \begin{align*}
           p_{\mu}= \kappa   \frac{ \det(\bar{\nabla}^2h+hI)}{\depsi_{t}(\cdot, h)|\bar{\nabla}h+h\iota|^n }. \end{align*}

  \section{A solution to the Musielak-Orlicz-Gauss image problem}\label{solution-8-25} \setcounter{equation}{0}

  Our goal in this section is to provide solutions to Problems \ref{MOGIP} and \ref{MOGIP-p}, mainly under the condition that $G$ is strictly decreasing on its second variable. Let  $\lambda\in\cM$ and  $\mu$ be nonzero finite Borel measures on $\sphere$.  Let $\wp:  \sphere \times (0, \infty)\to (0,\infty)$ be a  continuous function and $\Psi\in \cC_{d}$. Consider the following two  optimization problems: 
  \begin{align} 
       & \inf\left\{\widetilde{V}_{\depsi,\mu}(Q): \widetilde{V}_{G,\lambda}(Q^*)=\widetilde{V}_{G,\lambda}(B^n)\ \mathrm{and}\ Q\in
        \cKo\right\}, \label{geom} \\    &\inf\left\{\widetilde{V}_{\depsi,\mu}(f): \widetilde{V}_{G,\lambda}(\langle f\rangle^*)=\widetilde{V}_{G,\lambda}(B^n)\ \mathrm{and}\ f \in C^{+}\left(S^{n-1}\right)\right\},  \label{fun} \end{align} where $\widetilde{V}_{G,\lambda}(K)$ and
$\widetilde{V}_{G,\lambda}(f)$ are given in \eqref{vlam} and
\eqref{vlam2}, i.e., 
$$
  \VGL(K)=\int_{\sphere} \wp(\xi, \rho_K(\xi))\, d\lambda(\xi) \ \ \mathrm{and} \ \   \VGL(f)=\int_{S^{n-1}}G(\xi, f(\xi))\,d\lambda(\xi).$$ Recall that $\psi_{\xi}(t)=\depsi(\xi, t)$ and $\psi_{\xi}^{-1}$ is its inverse on $(0, \infty)$.

   The following lemma plays important roles in solving  Problem \ref{MOGIP-p}.  

\bl\label{m1} Let  $\lambda\in\cM$ and  $\mu$ be nonzero finite Borel measures on $\sphere$.  Suppose that  $\wp\in \cC$ and $\Psi\in \cC_{d}$ such that $C_{\Theta}(Q, \sphere)\neq 0$ for all $Q\in \cKo$.  If the optimization problem \eqref{geom} admits a solution, say $K\in \cKo$, 
then $K_0=K^*$ is a solution to Problem \ref{MOGIP-p}, namely, the following holds: 
    \begin{align}\label{solution-20-8-2}
    \frac{\mu}{|\mu|}=\frac{ C_{\Theta}(K_0 ,\cdot)}{ C_{\Theta}(K_0, \sphere)}.
    \end{align}
\el
 
\begin{proof}  For any  $f\in C^+(\sphere)$  such that
$\widetilde{V}_{G,\lambda}(\langle
f\rangle^*)=\widetilde{V}_{G,\lambda}(B^n)$, it follows from \eqref{hk} that  $\langle \rho_{\langle f\rangle}\rangle=\langle f\rangle$, which further implies   $\widetilde{V}_{G,\lambda}(\langle \rho_{\langle
f\rangle}\rangle^*)=\widetilde{V}_{G,\lambda}(\langle
f\rangle^*)=\widetilde{V}_{G,\lambda}(B^n)$ and 
$\widetilde{V}_{\depsi,\mu}(\rho_{\langle
f\rangle})=\widetilde{V}_{\depsi,\mu}(\langle
f\rangle)$. On the other hand, by $\Psi\in \cC_{d}$ and
\eqref{le1}, one has 
\begin{align*} \widetilde{V}_{\depsi,\mu}(f)=\int_{S^{n-1}}\Psi\left(\xi,  f(\xi)\right) \mathrm{d} \mu(\xi) \ge  \int_{S^{n-1}}\Psi\left(\xi, \rho_{\langle f\rangle}(\xi)\right) \mathrm{d} \mu(\xi)  = \widetilde{V}_{\depsi,\mu}\left(\rho_{\langle f\rangle}\right). \end{align*}  Hence if $K\in \cKo$ solves the optimization problem \eqref{geom},   then $\rho_K\in C^+(\sphere)$ solves the  optimization problem  \eqref{fun}.

  Let $g\in C(\sphere)$ be an arbitrary continuous function on $\sphere$. As $\rho_K\in C^+(\sphere)$, for sufficiently small
$\varepsilon_1,\varepsilon_2, \varepsilon$, it is meaningful to define  \begin{align*} 
        f_{\varepsilon_1+\varepsilon,\varepsilon_2}(\xi)=\psi_{\xi}^{-1}\left(\psi_{\xi}(f_{\varepsilon_1,\varepsilon_2}(\xi))+\varepsilon
        g(\xi)\right) \ \ \mathrm{and}\ \ 
        f_{\varepsilon_1,\varepsilon_2+\varepsilon}(\xi)=\psi_{\xi}^{-1}\left(\psi_{\xi}(f_{\varepsilon_1,\varepsilon_2}(\xi))+\varepsilon\right),
        \end{align*} where $f_{\varepsilon_1, \varepsilon_2}$ is given by 
        
        \begin{align}\label{feb172}
        f_{\varepsilon_1, \varepsilon_2}(\xi) &=\psi_{\xi}^{-1}\left(\psi_{\xi}(\rho_{K}(\xi))+\varepsilon_1 g(\xi)+\varepsilon_2\right). \end{align} A more convenient formula for $f_{\varepsilon_1, \varepsilon_2}$ given in \eqref{feb172} is  \begin{align}\label{feb172-1-2020-12}
        \depsi(\xi, f_{\varepsilon_1, \varepsilon_2}(\xi)) &= \depsi(\xi, \rho_{K}(\xi))+\varepsilon_1 g(\xi)+\varepsilon_2. \end{align}
It follows from \eqref{rr} (with $f=f_{\varepsilon_1, \varepsilon_2}$) that  
 \begin{align}\label{variation-11-27-1---1}
        \frac{\partial}{\partial\varepsilon_1}\widetilde{V}_{G,\lambda}(\langle f_{\varepsilon_1,
        \varepsilon_2}\rangle^*)&=
        \lim_{\varepsilon\rightarrow 0} \frac{\widetilde{V}_{G,\lambda}(\langle f_{\varepsilon_1+\varepsilon, \varepsilon_2}\rangle^*)-\widetilde{V}_{G,\lambda}(\langle f_{\varepsilon_1,\varepsilon_2}\rangle^*)}{\varepsilon} \nonumber \\
        &=
        -\int_{\sphere} g(\xi)\,dC_{\Theta}(\langle f_{\varepsilon_1,\varepsilon_2}\rangle^*,
        \xi).
        \end{align} Similarly, one can also have   \begin{align}\label{variation-11-27-1---2}
        \frac{\partial}{\partial\varepsilon_2}\widetilde{V}_{G,\lambda}(\langle f_{\varepsilon_1, \varepsilon_2}\rangle^*) =-\int_{\sphere} \,dC_{\Theta}(\langle f_{\varepsilon_1,\varepsilon_2}\rangle^*, \xi) =-C_{\Theta}(\langle f_{\varepsilon_1,\varepsilon_2}\rangle^*, \sphere)\neq 0.
        \end{align}
  Note that $f_{\varepsilon_1, \varepsilon_2}$ depends continuously on $\varepsilon_1$ and $\varepsilon_2$. Hence, part i) in Proposition \ref{prop 11-28--1} implies the weak convergence of  $C_{\Theta}(\langle f_{\varepsilon_1, \varepsilon_2}\rangle^*,\cdot)$ on $\varepsilon_1$ and $\varepsilon_2$, respectively. Together with (\ref{variation-11-27-1---1}) and (\ref{variation-11-27-1---2}), $\VGL(\langle f_{\varepsilon_1, \varepsilon_2}\rangle^*)$ has a gradient which has rank $1$ and is continuous on $\varepsilon_1$ and $\varepsilon_2$. In particular, $\VGL(\langle f_{\varepsilon_1, \varepsilon_2}\rangle^*)$  is continuously differentiable on $\varepsilon_1$ and $\varepsilon_2$. Therefore, the method of Lagrange multipliers can be applied to the optimization problem \eqref{fun} to get a constant $\kappa=\kappa(g)$ such that
        \begin{align}\label{lagrange method}
        \frac{\partial}{\partial \varepsilon_i}\left(\widetilde{V}_{\depsi,\mu}(f_{\varepsilon_1, \varepsilon_2})+\kappa \big(\widetilde{V}_{G,\lambda}(\langle f_{\varepsilon_1, \varepsilon_2}\rangle^*) - \widetilde{V}_{G,\lambda}(B^n)\big)\right)\Big|_{\varepsilon_1=\varepsilon_2=0}=0, \ \ \ i=1, 2.
        \end{align}
It follows from \eqref{bi-polar--12}, \eqref{relation}, \eqref{hk} and \eqref{feb172} that
        \begin{align*}
        \langle f_{0,0}\rangle^*=\langle \rho_{K}\rangle^*=\langle 1/h_{K^*}\rangle^*=[h_{K^*}]=K^*.
        \end{align*} Hence, $\widetilde{V}_{G,\lambda}(\langle f_{0,0}\rangle^*)=\widetilde{V}_{G,\lambda}(K^*)=\widetilde{V}_{G,\lambda}(B^n).$  Together with \eqref{vlam}, \eqref{feb172-1-2020-12}, (\ref{variation-11-27-1---1}), (\ref{variation-11-27-1---2}), and \eqref{lagrange method},  one can easily have 
        \begin{align}\label{identity-2020-1}
        \int_{S^{n-1}} g(\xi)\,d\mu(\xi) = \kappa(g) \int_{\sphere} g(\xi)\,dC_{\Theta}(K^*,
        \xi) \ \ \ \mathrm{and}\ \  \  |\mu|= \kappa(g)  C_{\Theta}(K^*, \sphere).
        \end{align} In particular, $\kappa=\kappa(g)$ is a constant independent of the choice of $g\in C(\sphere)$:
        \begin{align*}
        \kappa = \frac{|\mu|}{C_{\Theta}(K^*,\sphere)}.
        \end{align*} Thus, by \eqref{identity-2020-1},  the following formula holds for any $g\in C(\Omega)$:  \begin{align*}
        \int_{S^{n-1}} g(\xi)\,d\mu(\xi) =\frac{|\mu|}{C_{\Theta}(K^*,\sphere)} \ \int_{\sphere} g(\xi)\,dC_{\Theta}(K^*,
        \xi).
        \end{align*} This concludes that, by letting $K_0=K^*$,  \eqref{solution-20-8-2} holds on $\cB$.  
\end{proof}

\bl \label{Lemma-8-23} Let $G\in \cG_d$ and  $\mu$ be a nonzero finite Borel measure on $S^{n-1}$ that is not concentrated on any closed hemisphere. Assume that $\{K_i\}_{i\in \N}\subseteq \cKo$ is a sequence such that  \begin{align} \label{sup-8-25}\sup_{i\in \N}\int_{\sphere}G(\xi, \rho_{K^*_i}(\xi))\,d\mu(\xi)<+\infty. \end{align} Then, the sequence $\{K_i\}_{i\in \N}$ is uniformly bounded, namely, there exists a finite constant $R$ such that $K_i\subseteq R\ball$ for all $i\in \N$. 
 \el 
\begin{proof} For each $i\in \N$, let $R_i=\max_{v\in\sphere}\rho_{K_{i}}(v)$  and $v_i\in \sphere$ be such that $R_i=\rho_{K_{i}}(v_i)$. We now claim that $\sup_{i\in \N} R_{i}<+\infty$ by contradiction.  Assume not,   a subsequence  $\{R_{i_j}\}_{j\in \N}$ of $\{R_i\}_{i\in \N}$ can be obtained so that $v_{i_j}\rightarrow v_0\in \sphere$ (due to the compactness of $\sphere$) and $\lim_{j\to\infty}R_{i_j}=+\infty$.  For  $v\in\sphere$ and $\beta \in (0,1)$, let 
   \begin{align*} \Sigma_{\beta}(v)=\{u\in S^{n-1}:u\cdot v\ge \beta \}.  \end{align*}   
  
  As  $\mu$ is not concentrated on any closed hemisphere, a simple argument by the monotone convergence theorem implies the existence of $\beta_0\in (0, 1)$ such that   $\mu(\Sigma_{\beta_0}(v_0))>0$. For any $\xi \in \Sigma_{\beta_0}(v_0)$ and $i\in \N$, one has $\rho_{K_{i}}(v_i)v_i\in K_i$ and hence $$h_{K_i}(\xi)\ge \rho_{K_{i}}(v_i)(\xi\cdot v_i)=R_i(\xi\cdot v_i).  $$ The continuity of the function $u\mapsto \xi\cdot u$ and the facts that $\xi\cdot v_0\geq \beta_0$ and $v_{i_j}\to v_0$ yield the existence of $j_0\in \N$ such that for all $j\geq j_0$, $$h_{K_{i_j}}(\xi)\ge  R_{i_j}(\xi\cdot v_{i_j})\geq \frac{R_{i_j} \beta_0}{2}.  $$
 By \eqref{bi-polar--12},  \eqref{vlam},
\eqref{sup-8-25} and $\wp
\in\cG_d$, one gets, for all $j\geq j_0$, 
        \begin{align*}
         + \infty &>\int_{\sphere}G(\xi, \rho_{K^*_{i_j}}(\xi))\,d\mu(\xi)=   \int_{\sphere} \wp(\xi, h_{K_{i_j}}(\xi)^{-1})\,d\mu(\xi)  \ge \int_{\Sigma_{\beta_0}(v_0)} \wp\Big(\xi, \frac{2} {R_{i_j} \beta_0}\Big)\,d\mu(\xi).   \end{align*} 
    As $\lim_{t\rightarrow 0^+} G(\xi, t)=+ \infty$ for each $\xi\in \sphere$ and $\lim_{j\to\infty}R_{i_j}=+\infty$, Fatou's lemma yields that  \begin{align*}   + \infty> \liminf_{j \to \infty}  \int_{\Sigma_{\beta_0}(v_0)} \wp\Big(\xi, \frac{2} {R_{i_j} \beta_0}\Big)\,d\mu(\xi) \geq  \int_{\Sigma_{\beta_0}(v_0)}  \liminf_{j \to \infty}  \wp\Big(\xi, \frac{2} {R_{i_j} \beta_0}\Big)\,d\mu(\xi) = + \infty.          \end{align*} 
        This is a contradiction and hence the sequence $\{K_i\}_{i\in \N}$ is uniformly bounded. \end{proof} 
 
 \bl \label{Lemma-8-23-2} Let $\depsi\in \cC_d$ be such that  \begin{align}\label{growth-8-25} \lim_{t\rightarrow 0^+} \depsi(\xi, t)=+\infty \ \ \mathrm{for\ each}\ \xi\in \sphere. \end{align}  Let $\mu$ be a nonzero finite Borel measure on $S^{n-1}$ that is not concentrated on any closed hemisphere.  Assume that  the sequence $\{K_i\}_{i\in \N}\subseteq \cKo$ is uniformly bounded such that \begin{align} \label{max-8-23} \sup_{i\in \N} \int_{\sphere}\Psi \left(\xi, \rho_{K_i}(\xi)\right)\,d\mu(\xi)<+\infty. \end{align} Then, there exists a subsequence of  $\{K_i\}_{i\in \N}$ which converges to some $L\in \cKo$.  
 \el 
 \begin{proof} Let  $R$ be a finite constant such that $K_i\subseteq R\ball$ for all $i\in \N$.   Applying the Blaschke selection theorem to $\{K_i\}_{i\in \N}$,  a convex compact set $L\subseteq \Rn$  and a subsequence of $\{K_i\}_{i\in \N}$ can be found (which will still be denoted by $K_i$), such that $K_i\to L$ in the Hausdorff metric.

 We now claim that $L\in\cKo$. Assuming the contrary,  there exist $w_0\in S^{n-1}$ and $\beta_1>0$ such that $0=h_{L}(w_0)= \lim_{i\rightarrow \infty} h_{K_i}(w_0)$ and $\mu(\Sigma_{\beta_1}(w_0))>0$, where the latter one follows from the fact that $\mu$ is a nonzero finite Borel measure  not concentrated on any closed
hemisphere. From \eqref{vlam}, \eqref{max-8-23}, $\depsi\in \cC_d$, and  $K_i\subseteq R \ball$ for all $i\in \N$, one has 
\begin{align}
           + \infty&>
            \liminf_{i\rightarrow \infty} \int_{\sphere}\Psi \left(\xi, \rho_{K_i}(\xi)\right)\,d\mu(\xi) \nonumber \\ &\ge  \liminf_{i\rightarrow \infty} \int_{\Sigma_{\beta_1}(w_0)}\Psi \left(\xi, \rho_{K_i}(\xi) \right)\,d\mu(\xi)+\int_{\sphere\setminus \Sigma_{\beta_1}(w_0)}\Psi \left(\xi, R\right)\,d\mu(\xi) \nonumber \\&
            \ge  \liminf_{i\rightarrow \infty} \int_{\Sigma_{\beta_1}(w_0)}\Psi \left(\xi, \rho_{K_i}(\xi) \right)\,d\mu(\xi)
            +\mu(\sphere\setminus \Sigma_{\beta_1}(w_0)) \min_{u\in\sphere} \Psi \left(u, R \right).\label{con-83}
        \end{align} If $v\in
        \Sigma_{\beta_1}(w_0)$, then $\beta_1 \rho_{K_i}(v) \leq \rho_{K_i}(v)v\cdot w_0\leq h_{K_i}(w_0).$ Thus $\rho_{K_i}\rightarrow 0$ uniformly on $\Sigma_{\beta_1} (w_0)$ as $i\to\infty$, due to $ \lim_{i\rightarrow \infty} h_{K_i}(w_0)=0$.  This further yields, by \eqref{growth-8-25}, for any $\xi\in \sphere$,  $$ \liminf_{i\rightarrow \infty}\Psi \left(\xi, \rho_{K_i}(\xi) \right)=+\infty.$$
  A contradiction  can be obtained by \eqref{con-83}:  \begin{align*} +\infty &
           >  \liminf_{i\rightarrow \infty} \int_{\Sigma_{\beta_1}(w_0)}\Psi \left(\xi, \rho_{K_i}(\xi) \right)\,d\mu(\xi)
            +\mu(\sphere\setminus \Sigma_{\beta_1}(w_0)) \min_{u\in\sphere} \Psi \left(u, R \right)\\&
            \ge   \int_{\Sigma_{\beta_1}(w_0)}  \liminf_{i\rightarrow \infty}\Psi \left(\xi, \rho_{K_i}(\xi) \right)\,d\mu(\xi)+\mu(\sphere\setminus \Sigma_{\beta_1}(w_0)) \min_{u\in\sphere} \Psi \left(u, R \right)=+\infty,
        \end{align*} where in the second inequality, we have used the Fatou's lemma to the nonnegative functions $\Psi \left(\xi, \rho_{K_i}(\xi) \right)-\Psi \left(\xi, R \right)$ on $\Sigma_{\beta_1}(w_0)$. This concludes that  $L \in \cKo$ as desired.  \end{proof}

 Now let us prove the existence of a solution to Problem \ref{MOGIP-p}. 

\bt\label{rl}  Let $\lambda\in \cM$ and  $\mu$ be two nonzero finite Borel measures on $S^{n-1}$ that are not concentrated on any closed hemisphere. There exists a $K\in \cKo$  such that
    \begin{align}\label{msol}
    \frac{\mu}{|\mu|}=\frac{ C_{\Theta}(K,\cdot)}{ C_{\Theta}(K, \sphere)}
    \end{align} if either $\depsi\in \cC_d$ satisfies \eqref{growth-8-25} and  $\wp\in \cG_d$, or $G\in \cC_d$ satisfies \eqref{growth-8-25}  and   $\depsi\in \cG_d$.  
\et  
  \begin{proof}  Under the assumptions on $G$ and $\depsi$, it follows from  Proposition \ref{prop 11-28--1} iii) that $C_{\Theta}(Q, \sphere)\neq 0$ for all $Q\in \cKo$.  In view of Lemma \ref{m1}, one only needs to find an  $L \in \cKo$ which solves the optimization problem \eqref{geom}, i.e., $L$ must satisfy that 
  $\widetilde{V}_{G,\lambda}(L^*)=\widetilde{V}_{G,\lambda}(B^n)$ and $ \widetilde{V}_{\depsi,\mu}(L)=\alpha$ with 
 \begin{align}\label{optimization-convex-body-11-28-110}
       \alpha:=\inf\left\{\widetilde{V}_{\depsi,\mu}(Q):\  \widetilde{V}_{G,\lambda}(Q^*)=\widetilde{V}_{G,\lambda}(B^n)\ \mathrm{and}\ Q\in
        \cKo \right\}.
        \end{align} It is clear that the infimum is taking over a nonempty subset of $\cKo$ because $\ball$ satisfies the desired constraint condition in \eqref{optimization-convex-body-11-28-110}. In particular, this shows $$\alpha\leq \widetilde{V}_{\depsi,\mu}(\ball)=\int_{\sphere} \depsi(\xi, 1)\,d\mu(\xi)<+\infty.$$ Moreover, for each $\widetilde{Q}\in \cKo$, the   function $c\mapsto \widetilde{V}_{G,\lambda}(c\widetilde{Q}^*)$ is continuous on $c\in (0, \infty)$ and $ \widetilde{V}_{G,\lambda}(c_0\widetilde{Q}^*)= \widetilde{V}_{G,\lambda}(\ball)$ for some  $$c_0\in \Big[\min_{\xi\in \sphere} \rho^{-1}_{\widetilde{Q}^*}(\xi), \max_{\xi\in \sphere} \rho^{-1}_{\widetilde{Q}^*}(\xi)\Big].$$   The latter statement can be seen from the following argument:  the fact that $G$ is strictly decreasing on its second variable, and $\min_{\xi\in \sphere} \rho^{-1}_{\widetilde{Q}^*}(\xi)\cdot \widetilde{Q}^*\subseteq\ball \subseteq \max_{\xi\in \sphere} \rho^{-1}_{\widetilde{Q}^*}(\xi)\cdot \widetilde{Q}^*$ yield that
        $$ \widetilde{V}_{G,\lambda}\Big( \max_{\xi\in \sphere} \rho^{-1}_{\widetilde{Q}^*}(\xi)\cdot \widetilde{Q}^*\Big)\leq \widetilde{V}_{G,\lambda}(\ball)\leq \widetilde{V}_{G,\lambda}\Big( \min_{\xi\in \sphere} \rho^{-1}_{\widetilde{Q}^*}(\xi)\cdot \widetilde{Q}^*\Big). $$

In conclusion,  the optimization problem \eqref{optimization-convex-body-11-28-110} is well-defined and admits a  minimizing sequence, say $\left\{K_{i}\right\}_{i=1}^{\infty} \subseteq \cKo$, such that, by \eqref{vlam}, \begin{align}
        &\widetilde{V}_{G,\lambda}(B^n)=\widetilde{V}_{G,\lambda}(K_{i}^*)= \int_{\sphere}G(\xi, \rho_{K^*_i}(\xi))\,d\lambda(\xi)<+\infty, \label{sup-8-25--1}\\ 
        \label{lim}
& \alpha=\lim _{i \rightarrow \infty}
\widetilde{V}_{\depsi,\mu}\left(K_{i}\right)= \lim_{i\rightarrow \infty} \int_{\sphere}\Psi \left(\xi, \rho_{K_i}(\xi)\right)\,d\mu(\xi).
       \end{align} 
       
 For the case when $\depsi\in \cC_d$ satisfies \eqref{growth-8-25} and  $\wp\in \cG_d$,  one sees that \eqref{sup-8-25--1} verifies \eqref{sup-8-25}, and then Lemma \ref{Lemma-8-23} yields the uniform boundedness of  $\{K_i\}_{i\in \N}$. 
On the other hand, \eqref{lim} implies \eqref{max-8-23}, and Lemma \ref{Lemma-8-23-2} can be applied to obtain that (without loss of generality) $K_i\rightarrow L$ for some $L\in \cKo$. 
 For the case when  $G\in \cC_d$ satisfies \eqref{growth-8-25}  and   $\depsi\in \cG_d$, one sees that \eqref{lim}  verifies \eqref{sup-8-25}, and then Lemma \ref{Lemma-8-23} yields the uniform boundedness of $\{K_i^*\}_{i\in \N}$.  Similarly,  \eqref{sup-8-25--1} verifies \eqref{max-8-23}, and Lemma \ref{Lemma-8-23-2} can be applied to obtain that (without loss of generality) $K_i^*\rightarrow L^*$ for some $L\in \cKo$. In both cases, one has $K_i\to L\in \cKo$ and $K_i^*\to L^*$ as $i\rightarrow \infty$.   It is easily checked by the dominated convergence theorem and  $\widetilde{V}_{G,\lambda}(B^n)=\lim_{i\to \infty}\widetilde{V}_{G,\lambda}(K_{i}^*)$ for all $i\in \N$  that 
 \begin{align*} 
       & \alpha =\lim_{i\rightarrow \infty} \widetilde{V}_{\depsi,\mu}(K_i) =  \int_{S^{n-1}} \lim_{i\rightarrow
\infty}\Psi(\xi, \rho_{K_i}(\xi))\,d\mu(\xi) =
        \int_{S^{n-1}}\Psi(\xi, \rho_{L}(\xi))\,d\mu(\xi)=\widetilde{V}_{\depsi,\mu}(L),\\ & \widetilde{V}_{G,\lambda}(B^n) =  \int_{S^{n-1}} \lim_{i\rightarrow
\infty}G(\xi, \rho_{K^*_i}(\xi))\,d\mu(\xi) =
        \int_{S^{n-1}}G(\xi, \rho_{L^*}(\xi))\,d\mu(\xi)=\widetilde{V}_{G,\lambda}(L^*).
        \end{align*} Consequently $L\in \cKo$ solves the optimization problem \eqref{geom}.  By Lemma \ref{m1},  \eqref{msol} holds for $K=L^*$.    
        \end{proof}
        
The following corollary  is an easy consequence of Theorem \ref{rl} and Proposition \ref{prop 11-28--1} iii). 

\bc\label{rl-equiv}  Assume that either $\depsi\in \cC_d$ satisfies \eqref{growth-8-25} and  $\wp\in \cG_d$, or $G\in \cC_d$ satisfies \eqref{growth-8-25}  and   $\depsi\in \cG_d$. Let   $\lambda\in \cM$ be strictly positive on nonempty open subsets of $\sphere$ and $\mu$ be a  nonzero finite Borel measure on $\sphere$. Then the following statements are equivalent. 
\vskip 1mm \noindent i) The measure  $\mu$ on $S^{n-1}$ is not concentrated on any closed hemisphere. 
\vskip 1mm \noindent ii) There exists a $K\in \cKo$  such that
   $$
    \frac{\mu}{|\mu|}=\frac{ C_{\Theta}(K,\cdot)}{ C_{\Theta}(K, \sphere)}.
  $$
\ec

Recall that if  $\depsi\in\cG_I$, then
$\widetilde{\depsi}(\xi, t)=\depsi(\xi, \frac{1}{t})\in\cG_d$. Similarly, if  $\depsi\in\cC_I$, then
$\widetilde{\depsi}\in\cC_d$. Moreover, if $\depsi$ satisfies \begin{align}\label{growth-8-25-I} \lim_{t\rightarrow \infty} \depsi(\xi, t)=+\infty \ \ \mathrm{for\ each}\ \ \xi\in \sphere,\end{align} then $\widetilde{\depsi}$ satisfies \eqref{growth-8-25}.  Applying Theorem \ref{rl} and Corollary \ref{rl-equiv} to the triple $\widetilde{\Theta}=(G, \widetilde{\depsi}, \lambda)$, together with  \eqref{wc}, one can easily get a solution to the Musielak-Orlicz-Gauss image problem (i.e., Problem \ref{MOGIP}), which is summarized below.  

 \bt \label{solution-general-dual-Orlicz-main theorem-11-270}  Let $\lambda\in \cM$ and $\mu$ be two nonzero finite Borel measures on $S^{n-1}$ that are not concentrated on any closed hemisphere. Assume that  either $\depsi\in \cC_I$ satisfies \eqref{growth-8-25-I} and  $\wp\in \cG_d$, or $G\in \cC_d$ satisfies \eqref{growth-8-25}  and   $\depsi\in \cG_I$. Then, there exists a $K\in \cKo$  such that  
    \begin{align}\label{msol-8-5-1}
    \frac{\mu}{|\mu|}=\frac{\deV(K,\cdot)}{\deV(K, \sphere)}.
    \end{align} If, in addition, $\lambda$ is strictly positive on nonempty open subsets of $\sphere$, then  the assumption on $\mu$, i.e., $\mu$ is a nonzero finite Borel measure on $S^{n-1}$ that is not concentrated on any closed hemisphere, is also necessary for  \eqref{msol-8-5-1} holding true for some $K\in \cKo$.   
    \et
 
 Theorem \ref{solution-general-dual-Orlicz-main theorem-11-270} not only gives a Musielak-Orlicz generalization of \cite[Theorem 6.4]{GHWXY}, but also provides additional quite different assumptions on $G, \depsi$ (i.e., $G\in \cC_d$ satisfies \eqref{growth-8-25}  and   $\depsi\in \cG_I$) to guarantee the existence of solutions to the corresponding Minkowski type problems. In particular, the assumption that $G\in \cC_d$ satisfies \eqref{growth-8-25}  and   $\depsi\in \cG_I$ easily implies the existence of solutions to Problem \ref{MOMP-log-8-3}, 
 due to \eqref{equiv-8-25},  by letting $G=-\log t\in \cC_d$ which of course satisfies \eqref{growth-8-25}.  
\bt\label{sol-J-d}  Let $\lambda\in \cM$ and $\mu$ be two nonzero finite Borel measures on $S^{n-1}$ that are not concentrated on any closed hemisphere.  For $\Psi \in \cG_{I}$,  there exists a $K \in \cKo$  such that  
    \begin{align}
    \frac{\mu}{|\mu|}=\frac{\widetilde{J}_{\depsi, \lambda}(K, \cdot)}{\widetilde{J}_{\depsi, \lambda}\left(K,
S^{n-1}\right)}.\label{msol-8-25-1}
    \end{align}
     If, in addition, $\lambda$ is strictly positive on nonempty open subsets of $\sphere$, then  the assumption on $\mu$, i.e., $\mu$ is a nonzero finite Borel measure on $S^{n-1}$ that is not concentrated on any closed hemisphere, is also necessary for  \eqref{msol-8-25-1} holding true for some $K\in \cKo$.   \et

 The existence of  solutions to the Musielak-Orlicz-Aleksandrov problem (i.e., Problem \ref{MOMP-Alex}) is an easy consequence of  Theorem \ref{sol-J-d} by letting $\,d\lambda(\xi)=\,d\xi$.   
    \bc  Let $\Psi \in \cG_{I}$ and $\mu$ be a  nonzero finite Borel measure on $\sphere$.  The following two statements are equivalent. 
\vskip 1mm \noindent i) The measure  $\mu$  is not concentrated on any closed hemisphere.  \vskip 1mm \noindent ii) There exists a $K \in \cKo$  such that  $$
    \frac{\mu}{|\mu|}=\frac{\widetilde{J}_{\depsi}(K, \cdot)}{\widetilde{J}_{\depsi}\left(K,
S^{n-1}\right)}.$$ \ec

\section{A solution to the Musielak-Orlicz-Gauss image problem for even data}\label{even-8-25}  \setcounter{equation}{0} 
In this section, we will discuss the existence of solutions to the Musielak-Orlicz-Gauss image problem for even data. Most of the proofs in this sections follow along the lines similar to those in Section \ref{solution-8-25}, so we will mainly focus on the difference and modifications in the proofs.

Recall that a convex body $K\in \cKo$ is said to be origin-symmetric if $-x\in K$ for all $x\in K$. Denote by $\cKe\subseteq \cKo$ the collection of all origin-symmetric convex bodies. Let $C_e(\Omega)$ be the set of all even continuous functions defined on $\Omega\subseteq \sphere$, and $C^+_e(\Omega)$ contains all strictly positive functions in $C_e(\Omega)$.  Consider the following optimization problems:  \begin{align}   
        & \inf \left\{\widetilde{V}_{\depsi,\mu}(Q): \widetilde{V}_{G,\lambda}(Q^*)=\widetilde{V}_{G,\lambda}(\ball)\ \mathrm{and}\ Q\in
        \cKe\right\},   \label{geom-int-e-I} \\       &  \inf \left\{\widetilde{V}_{\depsi,\mu}(f): \widetilde{V}_{G,\lambda}(\langle f\rangle^*)=\widetilde{V}_{G,\lambda}(\ball) \ \mathrm{and}\ f \in C^+_e \left(S^{n-1}\right)\right\},  \label{fun-int-e-I} \\  &\alpha_s:= \sup \left\{\widetilde{V}_{\depsi,\mu}(Q): \widetilde{V}_{G,\lambda}(Q^*)=\widetilde{V}_{G,\lambda}(\ball)\ \mathrm{and}\ Q\in
        \cKe\right\},   \label{geom-int-e-I-825} \\       &  \sup \left\{\widetilde{V}_{\depsi,\mu}(f): \widetilde{V}_{G,\lambda}(\langle f\rangle^*)=\widetilde{V}_{G,\lambda}(\ball) \ \mathrm{and}\ f \in C^+_e \left(S^{n-1}\right)\right\}.  \label{fun-int-e-I-825} 
   \end{align}

\bl\label{m1-even} Let  $\lambda\in\cM$ and  $\mu$ be nonzero finite even Borel measures on $\sphere$.  Suppose that  $\wp\in \cC$ and $\Psi\in \cC$ such that $G(\xi, t)=G(-\xi, t)$ and $\depsi(\xi, t)=\depsi(-\xi, t)$ for all $(\xi, t)\in \sphere\times (0, \infty)$, and $C_{\Theta}(Q, \sphere)\neq 0$ for all $Q\in \cKe$.  
\vskip 1mm \noindent i) Let $\depsi\in \cC_d$.  If  \eqref{geom-int-e-I} admits a solution, say $K\in \cKe$, 
then $K_0=K^*\in \cKe$ is a solution to Problem \ref{MOGIP-p}, namely, the following holds: 
    \begin{align}\label{solution-20-8-2-e-I}
    \frac{\mu}{|\mu|}=\frac{ C_{\Theta}(K_0 ,\cdot)}{ C_{\Theta}(K_0, \sphere)}.
    \end{align}
 \noindent ii) Let $\depsi\in \cC_I$.  If  \eqref{geom-int-e-I-825} admits a solution, say $K\in \cKe$, 
then $K_0=K^*\in \cKe$ is a solution to Problem \ref{MOGIP-p}, namely, \eqref{solution-20-8-2-e-I} holds.  \el

\begin{proof}   
 Let $f\in C_e^+(\sphere)$  be such that
$\widetilde{V}_{G,\lambda}(\langle
f\rangle^*)=\widetilde{V}_{G,\lambda}(B^n)$. It follows from \eqref{hk} that  $\widetilde{V}_{G,\lambda}(\langle \rho_{\langle
f\rangle}\rangle^*)=\widetilde{V}_{G,\lambda}(B^n)$ and 
$\widetilde{V}_{\depsi,\mu}(\rho_{\langle
f\rangle})=\widetilde{V}_{\depsi,\mu}(\langle
f\rangle)$. If $\Psi\in \cC_{d}$, 
\eqref{le1} yields 
\begin{align*} \widetilde{V}_{\depsi,\mu}(f)=\int_{S^{n-1}}\Psi\left(\xi,  f(\xi)\right) \mathrm{d} \mu(\xi) \ge  \int_{S^{n-1}}\Psi\left(\xi, \rho_{\langle f\rangle}(\xi)\right) \mathrm{d} \mu(\xi)  = \widetilde{V}_{\depsi,\mu}\left(\rho_{\langle f\rangle}\right), \end{align*}  and similarly,  if $\Psi\in \cC_{I}$,  
\eqref{le1} yields  
\begin{align*} \widetilde{V}_{\depsi,\mu}(f)=\int_{S^{n-1}}\Psi\left(\xi,  f(\xi)\right) \mathrm{d} \mu(\xi) \le  \int_{S^{n-1}}\Psi\left(\xi, \rho_{\langle f\rangle}(\xi)\right) \mathrm{d} \mu(\xi)  = \widetilde{V}_{\depsi,\mu}\left(\rho_{\langle f\rangle}\right). \end{align*}  Hence,    $\rho_K\in C_e^+(\sphere)$ solves the  optimization problem  \eqref{fun-int-e-I}, if $K\in \cKe$ solves \eqref{geom-int-e-I}; while    
 $\rho_K\in C_e^+(\sphere)$ solves the  optimization problem  \eqref{fun-int-e-I-825}, if $K\in \cKe$ solves \eqref{geom-int-e-I-825}.

Let  $\wp\in \cC$ and $\Psi\in \cC_{I}\cup \cC_d$ satisfy the assumptions in Lemma \ref{m1-even}.  Let $g\in C_e(\sphere)$ be an arbitrary continuous function on $\sphere$. As $\rho_K\in C_e^+(\sphere)$, for sufficiently small
$\varepsilon_1$ and $\varepsilon_2$, it is meaningful to define $f_{\varepsilon_1, \varepsilon_2}$ as  in \eqref{feb172}.  
It follows from \eqref{rr} (with $f=f_{\varepsilon_1, \varepsilon_2}$) that  
 \begin{align} \label{variation-8-25-even}
        \frac{\partial}{\partial\varepsilon_1}\widetilde{V}_{G,\lambda}(\langle f_{\varepsilon_1,
        \varepsilon_2}\rangle^*)&=
        -\int_{\sphere} g(\xi)\,dC_{\Theta}(\langle f_{\varepsilon_1,\varepsilon_2}\rangle^*, 
        \xi), \\ \label{variation-even-8-25-2}
        \frac{\partial}{\partial\varepsilon_2}\widetilde{V}_{G,\lambda}(\langle f_{\varepsilon_1, \varepsilon_2}\rangle^*) &=-C_{\Theta}(\langle f_{\varepsilon_1,\varepsilon_2}\rangle^*, \sphere)\neq 0.
        \end{align} Again,  the method of Lagrange multipliers can be applied to the optimization problems \eqref{fun-int-e-I} or \eqref{fun-int-e-I-825} to get a constant $\kappa=\kappa(g)$, independent of $g$, such that
        \begin{align}\label{lagrange method-8-25-even}
        \frac{\partial}{\partial \varepsilon_i}\left(\widetilde{V}_{\depsi,\mu}(f_{\varepsilon_1, \varepsilon_2})+\kappa \big(\widetilde{V}_{G,\lambda}(\langle f_{\varepsilon_1, \varepsilon_2}\rangle^*) - \widetilde{V}_{G,\lambda}(B^n)\big)\right)\Big|_{\varepsilon_1=\varepsilon_2=0}=0, \ \ \ i=1, 2.
        \end{align} Note that $\widetilde{V}_{G,\lambda}(\langle f_{0,0}\rangle^*)=\widetilde{V}_{G,\lambda}(K^*)=\widetilde{V}_{G,\lambda}(B^n).$  Together with \eqref{vlam}, \eqref{feb172-1-2020-12}, (\ref{variation-8-25-even}), (\ref{variation-even-8-25-2}), and \eqref{lagrange method-8-25-even},  one can easily have, for all $g\in C_e(\sphere)$,  
          \begin{align*}
        \int_{S^{n-1}} g(\xi)\,d\mu(\xi) =\frac{|\mu|}{C_{\Theta}(K^*,\sphere)} \ \int_{\sphere} g(\xi)\,dC_{\Theta}(K^*,
        \xi).
        \end{align*} That is, \eqref{solution-20-8-2-e-I} holds by letting $K_0=K^*$.  
\end{proof} 
   
   The existence of solutions to Problem \ref{MOGIP-p} for even data is given in the following theorem.  
 
\bt\label{rl-even}  Let $\lambda\in \cM$ and  $\mu$ be two nonzero finite even Borel measures on $S^{n-1}$ that are not concentrated on any  great subsphere. Suppose that  $\wp\in \cC$ and $\Psi\in \cC$ such that $G(\xi, t)=G(-\xi, t)$ and $\depsi(\xi, t)=\depsi(-\xi, t)$ for all $(\xi, t)\in \sphere\times (0, \infty)$.

\vskip 2mm \noindent i)  If either $\depsi\in \cC_d$ satisfies \eqref{growth-8-25} and  $\wp\in \cG_d$, or $G\in \cC_d$ satisfies \eqref{growth-8-25}  and   $\depsi\in \cG_d$, then there exists a $K_0\in \cKe$  such that  \eqref{solution-20-8-2-e-I} holds.  
     
     If, in addition, $\lambda$ is strictly positive on nonempty open subsets of $\sphere$, then  the assumption on $\mu$, i.e., $\mu$ is a nonzero finite even Borel measure on $S^{n-1}$ that is not concentrated on any great subsphere, is also necessary for  \eqref{solution-20-8-2-e-I} holding true for some $K\in \cKe$.  
     
\vskip 2mm \noindent ii)  Assume that, in addition,  $\mu$ vanishes on great subspheres.  If $\depsi\in \cC_I$ and  $\wp\in \cG_d$, then there exists a $K_0\in \cKe$  such that  \eqref{solution-20-8-2-e-I} holds.

\vskip 2mm \noindent iii)  Assume that, in addition, $\mu$ vanishes on great subspheres and  there exists a constant $C\in (-\infty, \infty)$,  such that 
  \begin{align}\label{con-8-7-int} \inf_{v\in \sphere} \int_{S^{n-1}} \log \left|v \cdot
\xi\right| \,d\lambda(\xi) >C.\end{align}    If $\depsi\in \cC_I$, then there exists a $K\in \cKe$ such that  \begin{align}
    \frac{\mu}{|\mu|}=\frac{J_{\depsi, \lambda}(K, \cdot)}{J_{\depsi, \lambda}\left(K,
S^{n-1}\right)}.\label{con-8-7-int-21-1}
    \end{align}  
 \et  
\begin{proof}  For $v\in\sphere$ and $\beta\in (0,1)$, let
	$$
	\widehat{\Sigma}_{\beta}(v)=\{u\in S^{n-1}:| u\cdot v|\ge \beta\}. 
	$$
	
	\noindent  i) We first claim that the optimization problem \eqref{geom-int-e-I} has a solution under the assumptions $\depsi\in \cC_d$ satisfies \eqref{growth-8-25} and  $\wp\in \cG_d$. The case when $G\in \cC_d$ satisfies \eqref{growth-8-25}  and   $\depsi\in \cG_d$ follows along the same lines.

Following the proof of Theorem \ref{rl},  the optimization problem \eqref{geom-int-e-I} is well-defined and a limiting sequence $\{K_i\}_{i\in \N}$ can be found such that  $\widetilde{V}_{G, \lambda} (K_i^*)=\widetilde{V}_{G, \lambda} (\ball)$ for all $i\in \N$ and $\lim_{i\rightarrow \infty} \widetilde{V}_{\depsi, \mu}(K_i)\leq \widetilde{V}_{\depsi, \mu}(\ball).$  In particular, $G\in \cG_d$ and  $\{K_i\}_{i\in \N}\subseteq \cKe$ is a sequence satisfying \eqref{sup-8-25}. This shows that the sequence $\{K_i\}_{i\in \N}$ is uniformly bounded, which follows from an argument similar to the proof of Lemma \ref{Lemma-8-23}, mainly with $\cKo$ replaced by $\cKe$ and with the inner product replaced by its absolute value (consequently, with $\Sigma_{\beta}(\cdot)$  replaced by   $\widehat{\Sigma}_{\beta}(\cdot)$), respectively.  On the other hand, as $\depsi\in \cC_d$ satisfies  \eqref{growth-8-25} and $\{K_i\}_{i\in \N}\subseteq\cKe$ satisfies \eqref{max-8-23}, there exists a subsequence of  $\{K_i\}_{i\in \N}$ converging to some $L\in \cKe$; again this follows from an argument similar to the proof of Lemma  \ref{Lemma-8-23-2}. Without loss of generality, let $K_i\to L\in \cKe$ and then $K_i^*\to L^*$. Consequently, $\widetilde{V}_{G, \lambda} (L^*)=\widetilde{V}_{G, \lambda} (\ball)$ and $\lim_{i\rightarrow \infty} \widetilde{V}_{\depsi, \mu}(K_i)= \widetilde{V}_{\depsi, \mu}(L),$ namely, $L\in \cKe$ solves the optimization problem \eqref{geom-int-e-I}.  This, together with Lemma \ref{m1-even},  yields  $K_0=L^*\in \cKe$  satisfying  \eqref{solution-20-8-2-e-I}.

Recall that if  $\lambda\in \cM$ is strictly positive on nonempty open subsets of $\sphere$, then $C_{\Theta}(K, \cdot)$ is not concentrated on any closed hemisphere of $\sphere$. As $C_{\Theta}(K, \cdot)$ is an even measure, then  $C_{\Theta}(K, \cdot)$ is in fact not concentrated on any great subsphere. Consequently, if, in addition, $\lambda$ is strictly positive on nonempty open subsets of $\sphere$, then  the assumption on $\mu$, i.e.,  $\mu$ is a nonzero finite even Borel measure on $S^{n-1}$  not concentrated on any great subsphere,  is also necessary for  \eqref{solution-20-8-2-e-I} holding true for some $K\in \cKe$. 

\vskip 2mm \noindent ii) We first claim that the optimization problem \eqref{geom-int-e-I-825} has a solution under the assumptions $\depsi\in \cC_I$  and  $\wp\in \cG_d$.  Indeed, following the proof of Theorem \ref{rl}, the optimization problem \eqref{geom-int-e-I-825} is well-defined and a limiting sequence $\{K_i\}_{i\in \N}$ can be found such that  $\widetilde{V}_{G, \lambda} (K_i^*)=\widetilde{V}_{G, \lambda} (\ball)$ for all $i\in \N$ and  \begin{align} \alpha_s=\lim_{i\rightarrow \infty} \widetilde{V}_{\depsi, \mu}(K_i)=\lim _{i \rightarrow \infty} \int_{\sphere}\depsi \left(\xi, \rho_{K_i}(\xi)
       \right)\,d\mu(\xi) \geq \widetilde{V}_{\depsi, \mu}(\ball) >0.\label{lim-a-s}
\end{align} In particular, $G\in \cG_d$ and  $\{K_i\}_{i\in \N}\subseteq \cKe$ satisfy \eqref{sup-8-25}; this implies that the sequence $\{K_i\}_{i\in \N}$ is uniformly bounded (by an argument similar to the proof of Lemma \ref{Lemma-8-23}). Let $R$ be the constant such that $K_i\subseteq R\ball$ for all $i\in \N$.  The Blaschke selection theorem can be applied to get a compact convex set $L\subseteq\Rn$ and a subsequence of $\{K_i\}_{i\in \N}$ (which will still be denoted by $\{K_i\}_{i\in \N}$) such that $K_i\to L$ in the Hausdorff metric. 

Clearly $L$ is origin-symmetric.  If $L\notin \cKe$,  then there exists $w_0\in \sphere$, such that \begin{align}\label{L belong} L\subseteq w_{0}^{\perp}=\big\{x\in \R^n: \ x\cdot w_0=0\big\}.\end{align}  The fact that $\depsi\in \cG_I$ implies $0<\max_{u\in\sphere} \depsi \left(u, R\right):=C_1<\infty$. As $\mu$ is a nonzero finite even Borel measure  that vanishes on all great subspheres of
$S^{n-1}$, it can be checked that $$0=\mu( \sphere\cap w_0^{\perp})=\mu\Big(\bigcap\limits_{\beta\in (0, 1]}(\sphere\setminus
    \widehat{\Sigma}_{\beta}(w_0))\Big)=\lim_{\beta\rightarrow 0^+} \mu  (\sphere\setminus
    \widehat{\Sigma}_{\beta}(w_0)).$$  Let $\varepsilon>0$.  Then there exists  $\beta_{\varepsilon}\in (0, 1)$ such that   
\begin{align*}
 \mu(\sphere\setminus \widehat{\Sigma}_{ \beta_{\varepsilon}}(w_0))  
 <\frac{\varepsilon}{2C_1}.
        \end{align*}  As $\depsi\in \cG_I$ and $K_i\subseteq R\ball$ for all $i\in \N$, one has  \begin{align}
       \int_{\sphere\setminus \widehat{\Sigma}_{ \beta_{\varepsilon}}(w_0)}\Psi \left(\xi, \rho_{K_i}(\xi)
       \right)\,d\mu(\xi)\leq \int_{\sphere\setminus \widehat{\Sigma}_{ \beta_{\varepsilon}}(w_0)}\Psi \left(\xi, R
       \right)\,d\mu(\xi)<\frac{\varepsilon}{2}.\label{s-con-8-7}
        \end{align}
        
   It follows from \eqref{L belong} that $\lim_{i\rightarrow
\infty} h_{K_i}(w_0)=h_{L}(w_0)=0$.   This further implies that  $\rho_{K_i}\rightarrow 0$ uniformly on
$\widehat{\Sigma}_{\beta} (w_0)$ as $i\to\infty$ for any $\beta\in (0, 1)$ (see a similar argument in Lemma \ref{Lemma-8-23-2}). The dominated convergence theorem and $\depsi \in \cG_{I}$ yield the existence of $i_{\varepsilon} \in \N$, such that,  for
all $i\ge i_{\varepsilon} $,
\begin{align*}
       \int_{\widehat{\Sigma}_{\beta_{\varepsilon}}(w_0)}\depsi \left(\xi, \rho_{K_i}(\xi)
       \right)\,d\mu(\xi)<\frac{\varepsilon}{2}. 
        \end{align*} Together with \eqref{s-con-8-7}, one sees, for all $i\ge i_{\varepsilon} $,  
\begin{align*}
       \int_{\sphere}\depsi \left(\xi, \rho_{K_i}(\xi)
       \right)\,d\mu(\xi)<\varepsilon.
        \end{align*}  
 Taking \eqref{lim-a-s} into account, one gets a contradiction as follows: $$0= \alpha_s=\lim _{i \rightarrow \infty} \int_{\sphere}\depsi \left(\xi, \rho_{K_i}(\xi)
       \right)\,d\mu(\xi) =\int_{\sphere}\depsi \left(\xi, \rho_{L}(\xi)
       \right)\,d\mu(\xi) >0.$$ This concludes that $L\in \cKe$.

        In conclusion, one gets an origin-symmetric convex body $L\in \cKe$ such that $K_i\to L$ and then $K_i^*\to L^*$.    Moreover,   $\alpha_s= \widetilde{V}_{\depsi,\mu}\left(L\right)$ and $ \widetilde{V}_{G,\lambda}(L^*)=\widetilde{V}_{G,\lambda}(\ball)$, namely, $L\in \cKe$ solves  the optimization problem \eqref{geom-int-e-I-825}.  This, together with Lemma \ref{m1-even},  yields   $K_0=L^*\in \cKe$  satisfying  \eqref{solution-20-8-2-e-I}. 
        
    \vskip 2mm \noindent iii) In view of \eqref{pd-int},  to find a $K\in \cKe$ satisfying \eqref{con-8-7-int-21-1}, one needs to solve the optimization problem \eqref{geom-int-e-I-825} under the case $G=\log t$ (or, equivalently, $G=-\log t$ which can be clearly seen from \eqref{nn} and \eqref{pd-int}). In this case, $\widetilde{V}_{G,\lambda}(\cdot)$ has to be replaced by  $ \mathcal{E}_{\lambda}(K)$ defined in \eqref{rela}:  \begin{align*} \mathcal{E}_{\lambda}(K)=
   \widetilde{V}_{\log,\lambda}(K^*)=\int_{S^{n-1}} \log \rho_{K^*}(\xi) d \lambda(\xi).
    \end{align*} To be more precise,   the optimization problem \eqref{geom-int-e-I-825} now becomes 
    \begin{align}\label{max-21-1}
        \alpha_s:=\sup\left\{\widetilde{V}_{\widetilde{\depsi},\mu}(Q): \mathcal{E}_{\lambda}(Q)=0\
\mathrm{and}\ Q\in
        \cKe\right\}.
        \end{align} 
 Following the proof of Theorem \ref{rl}, the optimization problem \eqref{max-21-1} is well-defined and a limiting sequence $\{K_i\}_{i\in \N}$ can be found such that 
 $\mathcal{E}_{\lambda}(K_{i})=0$ for all $i\in \N$ and \begin{align*} \alpha_s=
\lim _{i \rightarrow \infty}
\widetilde{V}_{\widetilde{\depsi},\mu}\left(K_{i}\right)=\lim _{i \rightarrow \infty} \int_{\sphere}\widetilde{\Psi } \left(\xi, \rho_{K_i}(\xi)
       \right)\,d\mu(\xi)\geq \widetilde{V}_{\widetilde{\depsi},\mu}\left(\ball\right)>0.
\end{align*} 
 The sequence $\{K_i\}_{i\in \N}$ is uniformly bounded. To see this, let  $R_i=\max_{v\in\sphere}\rho_{K_{i}}(v)=\rho_{K_{i}}(v_i)$ and  $v_i\rightarrow v_0\in \sphere$. Assume that $\sup_{i\in \N} R_i=\infty$.  As $K_i\in \cKe$ for $i\in \N$, one has $h_{K_{i}}(\xi) \geq R_{i}\left|v_{i} \cdot \xi\right|$  for all $\xi \in S^{n-1}$.  By \eqref{bi-polar--12},  \eqref{rela}, \eqref{con-8-7-int}, and $ \mathcal{E}_{\lambda}\left(K_{i}\right)=0$ for all $i\in \N$, one gets, for all $i\in \N$, 
\begin{align*}
            0   =\int_{S^{n-1}} \log h_{K_{i}}(\xi)
d\lambda(\xi) \geq  \int_{S^{n-1}} \log \left|v_{i} \cdot
\xi\right| \,d\lambda(\xi)+\lambda(S^{n-1}) \log R_{i} \ge  C+\lambda(S^{n-1}) \log R_{i}. 
        \end{align*}
     Consequently,    a contradiction can be obtained as follows:   \begin{align*}+\infty=\sup_{i\in \N} \log R_{i} \leq \frac{- C} {\lambda(S^{n-1}) } <\infty. 
        \end{align*}  Hence $\sup_{i\in \N} R_i<\infty$ and the sequence $\{K_i\}_{i\in \N}$ is uniformly bounded. The rest of the proof is then identical to the proof for ii).   \end{proof}

  The existence of solutions to the Musielak-Orlicz-Gauss image problem (i.e., Problem \ref{MOGIP}) for even measures can be obtained by applying Theorem \ref{rl-even} to the triple $\widetilde{\Theta}=(G, \widetilde{\depsi}, \lambda)$ and  by \eqref{wc}. 
 
\bt\label{rl-even-ee}  Let $\lambda\in \cM$ and  $\mu$ be two nonzero finite even Borel measures on $S^{n-1}$ that are not concentrated on any  great subsphere. Suppose that  $\wp\in \cC$ and $\Psi\in \cC$ such that $G(\xi, t)=G(-\xi, t)$ and $\depsi(\xi, t)=\depsi(-\xi, t)$ for all $(\xi, t)\in \sphere\times (0, \infty)$. 

\vskip 2mm \noindent i) If either $\depsi\in \cC_I$ satisfies \eqref{growth-8-25-I} and  $\wp\in \cG_d$, or $G\in \cC_d$ satisfies \eqref{growth-8-25}  and   $\depsi\in \cG_I$,  there exists a $K_0\in \cKe$  such that    \begin{align}\label{solution-20-8-2-e-I-II}
    \frac{\mu}{|\mu|}=\frac{\deV(K_0 ,\cdot)}{\deV(K_0, \sphere)}.
    \end{align} 
   If, in addition, $\lambda$ is strictly positive on nonempty open subsets of $\sphere$, then  the assumption on $\mu$, i.e., $\mu$ is a nonzero finite even Borel measure on $S^{n-1}$ that is not concentrated on any great subsphere, is also necessary for  \eqref{solution-20-8-2-e-I-II} holding true for some $K_0\in \cKe$.  
   
   \vskip 2mm \noindent ii)  Assume that, in addition,  $\mu$ vanishes on great subspheres.  If $\depsi\in \cC_d$ and  $\wp\in \cG_d$, then there exists a $K_0\in \cKe$  satisfying  \eqref{solution-20-8-2-e-I-II}.
\et  

 The existence of solutions to Problem \ref{MOMP-log-8-3} for even measures can be established as well. Part i) of the following theorem is obtained from \eqref{equiv-8-25} and  Theorem \ref{rl-even-ee}, and  by letting $G=-\log t\in \cC_d$ which satisfies \eqref{growth-8-25}. Part ii) of the following theorem is obtained by \eqref{J-relation} and by  applying  Part iii) of Theorem \ref{rl-even} to the the function $\widetilde{\depsi}(\xi, t)=\depsi\big(\xi, 1/t\big)$ instead of $\depsi$ itself. 
 
\bt  \label{th-21-1} Let $\lambda\in \cM$ and $\mu$ be two nonzero finite even Borel measures on $S^{n-1}$ that are not concentrated on any great subsphere.
Let $\depsi\in \cC$  be such that $\depsi(\xi, t)=\depsi(-\xi, t)$ for all $(\xi, t)\in \sphere\times (0, \infty)$.

\vskip 2mm \noindent i)  If $\depsi\in \cG_I$, then there exists a $K\in \cKe$ such that \begin{align}
    \frac{\mu}{|\mu|}=\frac{\widetilde{J}_{\depsi, \lambda}(K, \cdot)}{\widetilde{J}_{\depsi, \lambda}\left(K,
S^{n-1}\right)}.\label{sol-int-e-8-7}
    \end{align} 
If, in addition, $\lambda$ is strictly positive on nonempty open subsets of $\sphere$, then the assumption on $\mu$, i.e., $\mu$ is not concentrated on any great subsphere of $\sphere$, is also necessary for  \eqref{sol-int-e-8-7} holding for some $K\in \cKe$.  

\vskip 2mm \noindent ii)  Assume that, in addition, $\mu$ vanishes on great subspheres and  there exists a constant $C\in (-\infty, \infty)$  such that \eqref{con-8-7-int} holds.   If $\depsi\in \cC_d$, then there exists a $K\in \cKe$ satisfying \eqref{sol-int-e-8-7}.
 \et

Note that $\int_{\{\xi:v\cdot\xi\neq0\}}\log |\xi\cdot v|\,d\xi=C$ is independent of $v$ and $C$ is finite.
Applying Theorem \ref{th-21-1} to the measure $\,d\lambda(\xi)=\,d\xi$, one can obtain a solution to the Musielak-Orlicz-Aleksandrov problem (i.e., Problem \ref{MOMP-Alex}) for even measures.   
   \bc \label{MMM-2021-1-16}  Let  $\mu$ be a nonzero finite even Borel measure on $S^{n-1}$ and  $\depsi\in \cG_I\cup \cG_d$  be such that $\depsi(\xi, t)=\depsi(-\xi, t)$ for all $(\xi, t)\in \sphere\times (0, \infty)$.

\vskip 1mm \noindent i) If $\depsi\in \cG_I$, then 
there exists a $K\in \cKe$ such that \begin{align}
    \frac{\mu}{|\mu|}=\frac{\widetilde{J}_{\depsi}(K, \cdot)}{\widetilde{J}_{\depsi}\left(K,
S^{n-1}\right)}\label{sol-int-e-8-7-sss}
    \end{align} if and only if   $\mu$ is a nonzero finite even Borel measure on $S^{n-1}$ that is not concentrated on any great subsphere of $\sphere$.

 \vskip 1mm \noindent ii)  Assume that, in addition,  $\mu$  vanishes on
any great subspheres of $S^{n-1}$.    If $\Psi \in \cG_d$, then there exists $K\in \cKe$ such that  \eqref{sol-int-e-8-7-sss} holds.  
\ec

\vskip 2mm \noindent  {\bf Acknowledgement.}  The research of QH has been supported by NSFC (No. 11701219) and AARMS postdoctoral fellowship. The research of DY has been supported by an NSERC grant, Canada. The research of  BZ has been supported by NSFC (No.\ 11971005).  The authors are greatly indebted to Dr. Shaoxiong Hou for the introduction of the Musielak-Orlicz functions.

\begin{small}
 
\end{small}

\vskip 2mm \noindent
Qingzhong Huang,   \ {\small \tt hqz376560571@163.com} \\
{\em School of Mathematics, Physics and Information Engineering, Jiaxing University, Jiaxing, Zhejiang, 314001,
China}

\vskip 2mm \noindent Sudan Xing, \ \ \ {\small \tt sxing@ualberta.ca}\\
{\em Department of Mathematical and Statistical Sciences, University of Alberta, Edmonton,  Alberta T6G 2G1, Canada
} 

\vskip 2mm \noindent Deping Ye, \ \ \ {\small \tt deping.ye@mun.ca}\\
{\em Department of Mathematics and Statistics, Memorial
University of Newfoundland, St. John's, Newfoundland A1C 5S7,
Canada
}

\vskip 2mm \noindent Baocheng Zhu, \ \ \ {\small \tt bczhu@snnu.edu.cn}\\
{ \em School of Mathematics and Statistics, Shaanxi Normal
University, Xi'an, 710062, China}

\end{document}